\documentclass[3p,,preprint,12pt]{elsarticle}

\usepackage{titlecaps}
\Addlcwords{the of into via for and of on in an to hp-finite with}

\usepackage{tabulary,xcolor}
\usepackage{amsfonts,amsmath,amssymb}
\usepackage[T1]{fontenc}
\makeatletter
\let\save@ps@pprintTitle\ps@pprintTitle
\def\ps@pprintTitle{\save@ps@pprintTitle\gdef\@oddfoot{\footnotesize\itshape \null\hfill\today}}
\def\hlinewd#1{%
  \noalign{\ifnum0=`}\fi\hrule \@height #1%
  \futurelet\reserved@a\@xhline}

\AtBeginDocument{\ifNAT@numbers \biboptions{sort&compress}\fi}
\makeatother

\usepackage{ifluatex}
\ifluatex
\usepackage{fontspec}
\defaultfontfeatures{Ligatures=TeX}
\usepackage[]{unicode-math}
\unimathsetup{math-style=TeX}
\else 
\usepackage[utf8]{inputenc}
\fi 
\ifluatex\else\usepackage{stmaryrd}\fi

\usepackage{url,multirow,morefloats,floatflt,cancel,tfrupee}
\makeatletter

\AtBeginDocument{\@ifpackageloaded{textcomp}{}{\usepackage{textcomp}}}
\makeatother
\usepackage{colortbl}
\usepackage{xcolor}
\usepackage{pifont}
\usepackage[nointegrals]{wasysym}
\makeatletter

\def\mcWidth#1{\csname TY@F#1\endcsname+\tabcolsep}

\def\cAlignHack{\rightskip\@flushglue\leftskip\@flushglue\parindent\z@\parfillskip\z@skip}
\def\rAlignHack{\rightskip\z@skip\leftskip\@flushglue \parindent\z@\parfillskip\z@skip}

\if@twocolumn\usepackage{dblfloatfix}\fi 
\AtBeginDocument{
\expandafter\ifx\csname eqalign\endcsname\relax
\def\eqalign#1{\null\vcenter{\def\\{\cr}\openup\jot\m@th
  \ialign{\strut$\displaystyle{##}$\hfil&$\displaystyle{{}##}$\hfil
      \crcr#1\crcr}}\,}
\fi
}

\let\lt=<
\let\gt=>
\def\processVert{\ifmmode|\else\textbar\fi}

\@ifundefined{subparagraph}{
\def\subparagraph{\@startsection{paragraph}{5}{2\parindent}{0ex plus 0.1ex minus 0.1ex}%
{0ex}{\normalfont\small\itshape}}%
}{}

\newcommand\role[1]{\unskip}
\newcommand\aucollab[1]{\unskip}
  
\@ifundefined{tsGraphicsScaleX}{\gdef\tsGraphicsScaleX{1}}{}
\@ifundefined{tsGraphicsScaleY}{\gdef\tsGraphicsScaleY{.9}}{}
\def\checkGraphicsWidth{\ifdim\Gin@nat@width>\linewidth
    \tsGraphicsScaleX\linewidth\else\Gin@nat@width\fi}

\def\checkGraphicsHeight{\ifdim\Gin@nat@height>.9\textheight
    \tsGraphicsScaleY\textheight\else\Gin@nat@height\fi}

\def\fixFloatSize#1{}
\let\ts@includegraphics\includegraphics

\def\inlinegraphic[#1]#2{{\edef\@tempa{#1}\edef\baseline@shift{\ifx\@tempa\@empty0\else#1\fi}\edef\tempZ{\the\numexpr(\numexpr(\baseline@shift*\f@size/100))}\protect\raisebox{\tempZ pt}{\ts@includegraphics{#2}}}}

\AtBeginDocument{\def\includegraphics{\@ifnextchar[{\ts@includegraphics}{\ts@includegraphics[width=\checkGraphicsWidth,height=\checkGraphicsHeight,keepaspectratio]}}}

\def\URL#1#2{\@ifundefined{href}{#2}{\href{#1}{#2}}}

\def\UrlOrds{\do\*\do\-\do\~\do\'\do\"\do\-}%
\g@addto@macro{\UrlBreaks}{\UrlOrds}
\makeatother

\usepackage{floatrow}

\usepackage{pgfplots,algorithmic,algorithm}
\pgfplotsset{compat=1.5}
\usepackage{booktabs}
\usepackage[frozencache]{minted}
\usepackage[toc,page]{appendix}

\usepackage[textsize=small]{todonotes}

\usepackage{pifont}

\usepackage{graphicx,amssymb,amsmath,amsthm}

\newcommand{\bx}{\mathbf{x}}
\newcommand{\by}{\mathbf{y}}
\newcommand{\bw}{\mathbf{w}}
\newcommand{\bW}{\mathbf{W}}
\newcommand{\bb}{\mathbf{b}}

\newcommand{\RR}{\mathbb{R}}

\theoremstyle{plain}
\newtheorem{assumption}{Assumption}
 
\theoremstyle{plain}
\newtheorem{remark}{Remark}

\usepackage{float}
\journal{Journal Name}

\usepackage{hyperref}
\usepackage[nameinlink]{cleveref}
\usepackage{booktabs}

\usepackage{subcaption}
\newtheorem{theorem}{Theorem}

\begin{document}

\begin{frontmatter}

\title{The Neural Network Approach to Inverse Problems in Differential Equations\tnoteref{t1}}

 \author[label1]{Kailai Xu}
 \ead{kailaix@stanford.edu}
 \address[label1]{Institute for Computational and Mathematical Engineering, Stanford University, Stanford, CA, 94305}
 \tnotetext[t1]{The first author thanks the Stanford Graduate Fellowship in Science \& Engineering and the 2018 Schlumberger Innovation Fellowship for the financial support.}
 
 \author[label2]{Eric Darve}
 \ead{darve@stanford.edu}
 \address[label2]{Mechanical Engineering, Stanford University, Stanford, CA, 94305}
 
\begin{abstract}

We proposed a framework for solving inverse problems in differential equations based on neural networks and automatic differentiation. Neural networks are used to approximate hidden fields. We analyze the source of errors in the framework and derive an error estimate for a model diffusion equation problem. Besides, we propose a way for sensitivity analysis, utilizing the automatic differentiation mechanism embedded in the framework. It frees people from the tedious and error-prone process of deriving the gradients. Numerical examples exhibit consistency with the convergence analysis and error saturation is noteworthily predicted. We also demonstrate the unique benefits neural networks offer at the same time: universal approximation ability, regularizing the solution, bypassing the curse of dimensionality and leveraging efficient computing frameworks.

\end{abstract}

\begin{keyword}
Neural Networks \sep Inverse Problems \sep Partial Differential Equations


\end{keyword}
\end{frontmatter}

\section{Introduction}

\subsection{Background}

Inverse problems are well motivated and challenging in many science and engineering applications. Given the model differential equations and initial/boundary conditions, we can solve the forward problem by numerical PDE schemes and give predictions of data. The inverse problems, on the other hand, start with the measured data and aim at estimating the parameters in the models. The inverse problems are encountered widely in seismic imaging~\cite{tian2014bayesian}, medical tomography~\cite{arridge1999optical}, plasma diagnosis~\cite{almleaky1989density} and so on. 

Recent development in machine learning and neural networks has made the neural network approach for solving inverse problems revive. In this article, we propose a framework of calibration of unknown functions in the model using neural networks. The use of neural networks has many benefits. The primary advantage of the proposed method is that the framework is widely applicable; theoretically, it can learn any continuous functions under mild assumptions, and potentially non-continuous functions. Besides, the neural networks also offer implicit regularization and can overcome the curse of dimensionalities. Thanks to the active involvement of deep learning researchers in both academia and industries, many inverse problems can be very easily implemented, and the performance is  usually optimized for general purpose. 

However, one concern that a researcher may have is whether the method has theoretical guarantee. In this article, we seriously consider the error estimate in the neural network approach. We prove that for the diffusion model problem, the error is bounded by $\mathcal{O}(\Delta t^2 + h^2)$, where $\Delta t$ and $h$ depends on the collected data, under certain assumptions on the optimization. We also demonstrate it through our numerical experiments, which show the consistency with our analysis. To our best knowledge, this is the first time this kind of analysis and numerical demonstration has been done in literature. 

The guiding principle is the same as that suggested in \cite{stuart2010inverse}: \textit{avoid discretization until the last possible moment}. The motivation is obvious: the discretization will introduce numerical errors~(in this article we will use the term \textit{consistency error} for this type of error in the inverse problem) which are undesired. Instead, we approximate the unknown functions directly using neural networks and evaluate it at the numerical discretization points. We find that the deferred discretization yields a more stable calibration.

Finally, the neural network approach provides a natural way for ``sensitivity analysis''. We derive a ``sensitive region'' for an interesting physical quantity under a specified neural network model. We have intentionally avoided the term ``uncertainty quantification'' since the sensitivity is optimization dependent; that is, it depends on the optimization procedure. The ``sensitive region'' describes to what extent the solution is credible if we are most interested in the accuracy of a particular physical quantity.
\subsection{Related Works}

There has been much work in the field of inverse problems in differential equations. Here we mention a few works. One family of approach is called \textit{sparsity regularization}~\cite{0266-5611-33-6-060301}. In 2004, Daubechies \textit{et al}~\cite{daubechies2004iterative} provided a first theoretical treatment on sparsity regularization for ill-posed inverse problems and established the convergence of the iterative soft-thresholding algorithm. Sparsity regularization as a paradigm for solving inverse problems has gained much popularity. For a detailed discussion on sparsity regularization in inverse problems, see \cite{scherzer2009variational,schuster2012tackling,schuster2012regularization,daubechies2016sparsity}. Another popular approach is the \textit{Bayesian approach}. This approach introduces regularization by means of the prior information and can handle nonlinearity and non-Gaussianity~\cite{Untitled92:online}. Besides, it is possible to formulate uncertainty quantification based on the Bayesian approach~\cite{dashti2016bayesian}. For detailed description, see \cite{dashti2016bayesian,bui2012gentle,stuart2010inverse,talay1996probabilistic,tian2014bayesian}. \textit{Adjoint methods} are also popular approaches, especially in seismic tomography~\cite{patella2005geophysical}. The adjoints method gives an efficient way to evaluate the gradients of an interesting function $g$ with respect to parameters in PDEs~(also called \textit{design parameters}, a.k.a. \textit{control variables} or \textit{decision parameters})~\cite{johnson2007notes}. The gradients are then fed into an optimization algorithm such as conjugate-gradient methods~\cite{shewchuk1994introduction}, and thereafter the best design parameters for maximizing~(or minimizing) $g$ are found. For more details, see \cite{hall1982sensitivity,jameson2003aerodynamic,fichtner2006adjoint,plessix2006review,giles2000introduction}.

Using neural networks for inverse problems has been in the literature for a while. Back to 30 years ago, Ibrahim~\cite{elshafiey1991neural} considered using a neural network for solving simple inverse problems such as the Fredholm integral equations. Another paper~\cite{ogawa1998neural} published in 1998 considered the network inversion method for the Fredholm integral equations of the first kind. Unfortunately, due to the limited computation capacity, the authors only considered one- or two-layer neural networks, which have limited representation ability and therefore do not suffice for more complicated tasks. Recently, there has been great process in using neural networks for engineering problems~\cite{raissi2017physics,raissi2017physics2,pang2018neural,zhang2018quantifying,wang2018deepmd,han2018mean,zhang16reinforced,wu2018understanding}. Our work distinguishes from other works in that 
\begin{itemize}
    \item We abide by the ``\textit{avoid discretization until the last possible moment}'' principle. For example, for approximating an unknown multivariate scalar function, the neural network exactly takes in a position vector and outputs a number.
    \item Besides, the algorithm is designed to incorporate the physics, utilizing the existing numerical schemes instead of creating special schemes. This also enables ``\textit{hot-plugging}'', where users are still able to use sophisticated numerical schemes.
    \item  Another noteworthy feature is that algorithmically we can pretend we know the hidden fields and solve the forward problem. The objective function obtained is by comparing the obtained solution and observations. We can treat the algorithm as a blackbox: once fed with observations, we get the approximated unknown function immediately. This feature is called \textit{end-to-end}~\cite{2017arXiv170408305G} in the machine learning community.
    \item  Finally, under this framework, we can give convergence bound and classify the sources of error in terms of optimization, discretization, and observation. Pioneer works in \cite{pang2018fpinns} inspire the latter view.
\end{itemize}

We summarize the contribution of the article here
\begin{itemize}
    \item Propose a neural network framework for solving inverse problems. 
    \item Conduct error analysis for the framework. Particularly, we proved the convergence rate of the approach as data is collected in finer granularity. 
    \item Propose a way for sensitive analysis under the framework. 
    \item Solve several calibrating problems for linear and nonlinear differential equations. We also give general comments on the numerical results. 
\end{itemize}

\section{Calibrating Unknown Functions using Neural Networks}\label{sect:algo}

In this paper, we consider the evolution problems. Consider a physical system that is governed by a set of partial differential equations
\begin{equation}\label{equ:pde}
    u_t = \mathcal{L}(\bx, u; f)
\end{equation}
here $f$ is a unknown function which maps $\bx$ to a scalar or a vector. We have a set of observations~(possibly noisy) for $u(\bx, t)$. In this paper, we consider sequential snapshots of the function values, which can be viewed as samples from the set of distributions 
\begin{equation}
    \mathcal{M} = \{ (u(\bx, t)+\bw_0, u(\bx, t+\Delta t )+\bw_1, \ldots, u(\bx, t+(m-1)\Delta t)+\bw_{m-1}) | t\in\mathcal{T} \}
\end{equation}
where $\mathbf{w}_i$'s are random noise~(not necessarily i.i.d.) and $\mathcal{T}$ is a finite set of positive values~(sampling times). $\mathcal{M}$ can be seen as that we have taken a continuous sequence of snapshots at each time $t\in \mathcal{T}$. For later discussion, we consider a sample $M$ subject to $\mathcal{M}$ 
\begin{align}
    M &= \{ (U_0, U_1, \ldots, U_{m-1}) | t\in\mathcal{T} \} \\
     &= \{ (u(\bx, t)+\bW_0, u(\bx, t+\Delta t )+\bW_1, \ldots, u(\bx, t+(m-1)\Delta t)+\bW_{m-1}) |\\
     &\qquad\qquad\qquad t\in\mathcal{T}, \bx = \{x_1, x_2, \ldots, x_n\} \}
\end{align}

Our task is to find an appropriate $f_\theta$ such that the data are consistent with the model. Here $f_\theta$ is parametrized by $\theta$. The basic idea is first to propose a function $f_\theta$, and then do the forward simulation and see if the predicted result is consistent with the data. If not, we propose another $f_{\theta'}$. The transition of  $f_{\theta}$ to $f_{\theta'}$ is carried out by an first-order optimizer that uses automatic differentiation. The utilization of the gradient information will boost the performance tremendously compared to other zeroth order methods, which are intrinsically trial-and-error methods. 

Another key component is that we combine the robust and highly accurate numerical methods in partial differential equations with the neural networks. Assume that the numerical scheme for \cref{equ:pde} is 
\begin{equation}\label{equ:discretize}
    L(v^{N+m-1}, v^{N+m-2}, \ldots, v^{N}; f_{\theta}) = 0
\end{equation}
where $v^{N}(\bx) = v(\bx, t_N)$ is the numerical solution for $u(\bx, t_N)$ at time $t=t_N$. The scheme can be explicit, implicit or semi-implicit; it can also be linear or nonlinear. 

The loss function to minimize is given by 
\begin{equation}\label{equ:loss}
    \mathrm{loss} = \sum_{t_N\in \mathcal{T}} L(U^{N+m-1}, U^{N+m-2}, \ldots, U^{N}; f_{\theta})^2
\end{equation}
Note even if we minimize the loss function to zero, we are not guaranteed that $\tilde f = f$ since in general 
\begin{equation}\label{equ:tf}
    L(u^{N+m-1}, u^{N+m-2}, \ldots, u^{N}; f)\neq 0
\end{equation}
due to discretization error. However, if the error is small enough, the optimization drives \cref{equ:loss} to near zero, the noise is small enough, and the neural network approximation ability is sufficient, we expect $f_{\theta} \approx f$. In the next section, we will conduct a quantified analysis of those errors. 

We have to point out it is imperative to adopt the \textbf{deferred discretization} here. That is, we do not first learn $\{f_{\theta}(\bx_i)\}$ at the discretization points $\{\bx_i\}$ and then try to fit these values using a neural network. Instead, we first form a neural network and evaluate the neural network at the discretization points. The neural network structure will be 
\begin{equation}\label{equ:structure}
    \bx\in\RR^d \rightarrow \textrm{multilayer dense neural networks} \rightarrow y\in \RR
\end{equation}

Like most other inverse problems, such problems are usually ill-posed, and most of the time the stability is questionable~(we use the term \textit{ill-conditioned} here). That is, the learned function $\tilde f$ can be very sensitive to the data we have collected, especially those with noise. Neural networks offer a regularization effect in the deferred discretization setting. Indeed, if we use smooth (nonlinear) activation functions in the ``multilayer dense neural networks'' part in \cref{equ:structure}, the function $f_{\theta}:\RR^d\rightarrow \RR$ will be smooth. And if we appropriately regularize the neural network, the gradients of the neural networks can also be bounded~(for example, by adopting a projected gradient descent during the training process). That results in a ``well-behaved'' function. In this setting, if the unknown function is also well-behaved, the neural network will offer a good solution in general. The power of the neural network is actually when the unknown function is ill-behaved. In this case, the neural network will either fit the function well or signal us through the unusual behavior of the optimization or the extreme values in the weights. 

\section{Error Analysis and Sensitivity Analysis}

\subsection{Approximation of Unknown Functions by Neural Networks}

\paragraph{Neural Network Architecture} We adopt the standard dense neural networks, which can be expressed as a series of function compositions
\begin{equation}\label{equ:nn}
  \begin{aligned}
    \by_2(\bx) &= \tanh(\bW_1 \bx + \bb_1) \\
    \by_3(\by_2) &= \tanh(\bW_2 \by_2 + \bb_2) \\
    \ldots \\
    \by_{n_l}(\by_{n_l-1}) &= \tanh(\bW_{n_l-1} \by_{n_l-1} + \bb_{n_l-1})\\
    \by_{n_l+1}(\by_{n_l}) &= \bW_{n_l} \by_{n_l} + \bb_{n_l}\\
    f_\theta(\bx) &= \by_{n_l+1}(\by_{n_l}(\ldots (\by_2(\bx) ) ))
\end{aligned}
\end{equation}
here $\theta$ ensembles all the weights~(including biases) parameters 
\begin{equation}
    \theta = \{ \bW_1, \bW_2, \ldots,\bW_{n_l}, \bb_1, \bb_2, \ldots, \bb_{n_l}  \}
\end{equation}

We have used the activation functions \texttt{tanh} in this work, but other choices are also possible, such as ReLU, sigmoid, leaky ReLU and so on. Different neural network architectures are also possible, such as those that contain the convolutional layers~\cite{krizhevsky2012imagenet}, sparse convolutional neural networks~\cite{jia2018modified}, pooling layers, residual connections~\cite{he2016deep}, recurrent cells~\cite{hochreiter1997long} and so on. 

\paragraph{Uniformly bounded weights} Regularization techniques are extensively studied to combat overfitting in neural network researches. It was shown that the magnitude of the weights directly impacts the generalization ability of the neural networks. Intuitively, functions are simpler when they vary at a slower rate and thus generalize better~\cite{gouk2018regularisation}. This idea is explored in various articles and shown to be effective to achieve the state-of-the-art performance. Therefore, following the line of research  in machine learning, we constrain our neural networks to be equipped with uniformly bounded weights, i.e., there exists a positive constant $C>0$, such that for all training scenarios
\begin{equation}\label{equ:constraint}
    \|\bW_i\|_2 \leq C\qquad i = 1,2,3,\ldots,n_l
\end{equation}
Here we have chosen $p=2$ norm for $\bW_i$, but other choices such as $p=1$, $\infty$ are possible. In addition, we do not need to make any assumptions on the bias terms $\bb_i$. For technical reasons, we make the following assumption, which imposes a bound on the last bias term $\bb_{n_l}$

\begin{assumption}[Uniform Bound on Parameters]\label{assum:bound}
    There exists a positive constant $C>0$ such that 
    \begin{align}
        \|\bW_i\|_2 \leq & C\qquad i = 1,2,3,\ldots,n_l\\
        \|\bb_{n_l}\|_2 \leq & C 
    \end{align}
\end{assumption}

The constraints can be easily enforced through a projection step in the optimization phase and perform a variant of the projected (stochastic) gradient descent method. To be specific, we define a projection function
\begin{equation}
    \pi(\bW, C) = \frac{1}{\max\{1, \frac{\|\bW\|_2}{C} \}} \bW
\end{equation}
which will project the weights back to the closest matrix with bounded 2-norm $C$. Under this assumption, we are able to derive an upper bound on the first and second order~(and theoretically all orders, albeit the upper bound will go to infinity as the order increases) derivatives of the neural network approximation function.
\begin{theorem}\label{thm:uniform_bounded}
    Assume \cref{equ:constraint} holds. Then the neural network given by \cref{equ:nn} satisfies
    \begin{align}
        \left\|\frac{\partial f_\theta(\bx)}{\partial\bx} \right\|_2\leq & C^{n_l}\\
        \left\|\frac{\partial^2 f_\theta(\bx)}{\partial\bx^2} \right\|_2\leq & 2C^{n_l+1}\frac{1-C^{n_l-1}}{1-C}
    \end{align}
\end{theorem}

\begin{proof}
    We use the numerator layout for the matrix calculus. Using the fact that $\frac{\partial \tanh(x)}{\partial x} = 1-\tanh(x)^2$ we have according to the chain rule
    \begin{align}
        \frac{\partial f_\theta(\bx)}{\partial\bx}  & =\frac{{\partial {{\mathbf{y}}_{{n_l} + 1}}}}{{\partial {{\mathbf{y}}_{{n_l}}}}}\frac{{\partial {{\mathbf{y}}_{{n_l}}}}}{{\partial {{\mathbf{y}}_{{n_l} - 1}}}} \ldots \frac{{\partial {{\mathbf{y}}_2}}}{{\partial {\mathbf{x}}}} \\
        &= {{\mathbf{W}}_{{n_l}}}{{\mathbf{W}}_{{n_l} - 1}} \ldots {{\mathbf{W}}_1}(1 - {\mathbf{y}}_{{n_l}}^{} \otimes {\mathbf{y}}_{{n_l}}^{})(1 - {\mathbf{y}}_{{n_l} - 1}^{} \otimes {\mathbf{y}}_{{n_l} - 1}^{}) \ldots (1 - {\mathbf{y}}_2^{} \otimes {\mathbf{y}}_2^{})
    \end{align}
    here $\otimes$ denotes element-wise multiplication. Note that $\tanh(x)\in (-1,1)$, we immediately obtain that 
    \begin{equation}
        \left\|\frac{\partial f_\theta(\bx)}{\partial\bx} \right\|_2\leq {\left\| {{{\mathbf{W}}_{{n_l}}}} \right\|_2}{\left\| {{{\mathbf{W}}_{{n_l} - 1}}} \right\|_2} \ldots {\left\| {{{\mathbf{W}}_1}} \right\|_2} \leqslant {C^{{n_l}}}
    \end{equation}
    
    Likewise, we have 
    \begin{equation}\label{equ:p1}
  \begin{aligned}
        \frac{\partial^2 f_\theta(\bx)}{\partial\bx^2}  &= \frac{\partial }{{\partial {\mathbf{x}}}}\left( {\frac{{\partial {{\mathbf{y}}_{{n_l} + 1}}}}{{\partial {{\mathbf{y}}_{{n_l}}}}}\frac{{\partial {{\mathbf{y}}_{{n_l}}}}}{{\partial {{\mathbf{y}}_{{n_l} - 1}}}} \ldots \frac{{\partial {{\mathbf{y}}_2}}}{{\partial {\mathbf{x}}}}} \right)\\
        &=\sum\limits_{i = 2}^{{n_l}} {{{\mathbf{W}}_{{n_l}}}{{\mathbf{W}}_{{n_l} - 1}} \ldots {{\mathbf{W}}_1}\left( { - 2{\mathbf{y}}_i^{}\frac{{\partial {{\mathbf{y}}_i}}}{{\partial {\mathbf{x}}}}} \right)} \prod\limits_{j \ne i,j = 2}^{{n_l}} {} (1 - {\mathbf{y}}_j^{} \otimes {\mathbf{y}}_j^{})
    \end{aligned}
\end{equation}

    Since we have
    \[\frac{{\partial {{\mathbf{y}}_i}}}{{\partial {\mathbf{x}}}} = \frac{{\partial {{\mathbf{y}}_i}}}{{\partial {{\mathbf{y}}_{i - 1}}}}\frac{{\partial {{\mathbf{y}}_{i - 1}}}}{{\partial {{\mathbf{y}}_{i - 2}}}} \ldots \frac{{\partial {{\mathbf{y}}_2}}}{{\partial {\mathbf{x}}}} = {{\mathbf{W}}_{i - 1}}{{\mathbf{W}}_{i - 2}} \ldots {{\mathbf{W}}_1}\prod\limits_{j = 2}^i {} (1 - {\mathbf{y}}_j^{} \otimes {\mathbf{y}}_j^{})\]
    and therefore
    \[{\left\| {\frac{{\partial {{\mathbf{y}}_i}}}{{\partial {\mathbf{x}}}}} \right\|_2} \leqslant {C^{i - 1}}\]
    
 Insert the equation into \cref{equ:p1} we have
 \begin{equation}
     \left\|\frac{\partial^2 f_\theta(\bx)}{\partial\bx^2} \right\|_2 \leq 2{C^{{n_l}}}\sum\limits_{i = 2}^{{n_l}} {{C^{i - 1}}}  = 2{C^{{n_l} + 1}}\frac{{{C^{{n_l} - 1}} - 1}}{{C - 1}}
 \end{equation}
\end{proof}

\subsection{Approximation Theory of Neural Networks}

The approximation degree for one layer neural network is well understood in literature. One such result is~\cite{mhaskar1996neural}
\begin{theorem}\label{thm:1d}
     Let $1\leq d'\leq d$, $r\geq 1$, $d\geq 1$ be integers, $1\leq p\leq \infty$, $\varphi:\RR^{d'}\rightarrow\RR$ be infinitely many times continuously differentiable in some open sphere in $\RR^{d'}$. We further assume that there exists $\mathbf{b}$ in this sphere such that 
     \begin{equation}
         \mathbf{D}^\mathbf{k} \varphi(\mathbf{b})\neq 0,\quad \mathbf{k}\in \mathbb{Z}^{d'}, \mathbf{k}\geq \mathbf{0}
     \end{equation}
     Then there exist $d'\times d$ matrices $\{A_j\}_1^n$ with the following property. For any $f\in W^p_{r,d}$, there exist coefficients $a_j(f)$ such that 
     \begin{equation}
         \|f-\sum_{j=1}^n a_j(f) \varphi(A_j(\cdot)+\mathbf{b})\|_p \leq cn^{-r/d}\|f\|_{W^p_{r,d}}
     \end{equation}
 \end{theorem}
 Here
 \begin{equation}
     \|f\|_{W^p_{r,d}} = \sum_{0\leq |\mathbf{k}|\leq r} \|D^{\mathbf{k}}f\|_{p}
 \end{equation}
 with multi-integer $\mathbf{k} = (k_1, k_2, \ldots, k_s)$
 
 For example, if $d'=1$, and $\varphi=\sigma$ is the sigmoid function, then $\varphi$ satisfies the condition in the theorem. If $f\in W_{1, d}^p$~(and therefore $r=1$), we obtain a neural network with one hidden layer~(hidden size equals $n$). To achieve a predetermined error $\varepsilon>0$, we need
 \begin{equation}\label{equ:exp}
     n^{-1/d} = \mathcal{O}(\varepsilon) \Rightarrow n = \mathcal{O}(\varepsilon^{-d})
 \end{equation}
 This indicates that one layer neural network suffers from the curse of dimensionalities -- the number of required neurons in the hidden layer grows exponentially with the input dimension $d$.
 
One possible way to overcome the difficulty is to consider multilayer neural network. We will now consider this possibility. 
 
In the late 1950s, Kolmogorov proved in a series of papers the Kolmogorov Superposition Theorem that answers~(in the negative) Hilbert's 13th problem~\cite{liu2015kolmogorov}. It is a theorem about representing instead of approximating and quite deep. It says that for a continuous function defined on $[0,1]^d$, given appropriate activation function, a neural network with hidden sizes $\mathcal{O}(d)$ will be able to represent the function \textit{exactly}.
\begin{theorem}[Kolmogorov Superposition Theorem~\cite{lorentz1996constructive}]\label{thm:kol}
    There exist $d$ constants $\lambda_j>0$, $j=1$, $2$, $\ldots$, $d$, $\sum_{j=1}^d \lambda_j\leq 1$ and $2n+1$ strictly increasing continuous function $\phi_i$, $i=1$, $2$, $\ldots$, $2d+1$, which map $[0,1]$ to itself, such that every continuous function $f$ of $d$ variables on $[0,1]^d$ can be represented in the form 
    \begin{equation}\label{equ:super}
        f(x_1, \ldots, x_d) = \sum_{i=1}^{2d+1} g\left( \sum_{j=1}^d \lambda_j \phi_i(x_j) \right)
    \end{equation}
    for some $g\in C[0,1]$ depending on $f$.
\end{theorem}

Based on the above result, Maiorov and Pinkus~\cite{maiorov1999lower} proved that if the fixed number of units in the hidden layers are $6d+3$ and $3d$, the two-layer neural networks can approximate any function to arbitrary precision.
\begin{theorem}\label{thm:mai}
    There exists an activation function $\sigma$ which is analytical, strictly increasing and sigmoidal and has the following property: for any $f\in C[0,1]^d$ and $\varepsilon>0$, there exist constants $d_i$, $c_{ij}$, $\theta_{ij}$, $\gamma_i$ and vectors $\bw^{ij}\in \RR^d$ for which 
    \[\left| {f({\mathbf{x}}) - \sum\limits_{i = 1}^{6d + 3} {{d_i}\sigma \left( {\sum\limits_{j = 1}^{3d} {{c_{ij}}\sigma ({{\mathbf{w}}^{ij}} \cdot {\mathbf{x}} - {\theta _{ij}})}  - {\gamma _i})} \right)} } \right| < \varepsilon \]
    for all $\mathbf{x}=(x_1, \ldots, x_d)\in [0,1]^d$.
\end{theorem} 

\Cref{thm:kol,thm:mai} only implies that the approximation ability of two-layer neural network for appropriate activation functions. These activation functions can be quite pathological to achieve the desired accuracy. It is unclear if the activation functions used in practice, such as sigmoid functions, suffice for accurate approximation.  

For any $f\in C[0,1]^d$, let the corresponding superposition decomposition be \cref{equ:super}. Consider the family of the functions $f$ where the corresponding $g\in W_{1,1}^p$ and $\phi_i(x)\in W_{1,1}^p$, then by \cref{thm:1d}, there exist constants ${{w_{ri}}}$, ${{a_{ri}}}$, ${{b_{ri}}}$ and a constant $C>0$ independent of $r$, such that

\begin{equation}\label{equ:gt}
\left| {g(t) - \sum\limits_{i = 1}^r {{w_{ri}}\sigma ({a_{ri}}t + {b_{ri}})} } \right| \leqslant \frac{{C||g|{|_{W_{1,1}^p}}}}{r}
\end{equation}

We call the function family of $f$ a \textit{regular family} if $w_{ri}$, $a_{ri}$ are bounded independent of $r$, i.e., there exists a constant $\tilde C$, s.t.
\begin{equation}
    |w_{ri}|\leq \tilde C, |a_{ri}|\leq \tilde C\qquad i = 1,2,\ldots, r, \qquad \forall r 
\end{equation}

\begin{remark}
    The assumptions that $g\in W_{1,1}^p$ and $\phi_i(x)\in W_{1,1}^p$ are critical for the error bound. Note that even strictly increasing continuous functions can be very pathological. An example is the the Cantor function $\mathcal{D}(x)$. $\mathcal{D}(x)$ is a nondecreasing continuous function but does not have a weak derivative. $x+\mathcal{D}(x)$ will then be a strictly increasing continuous function but has no weak derivatives. We want to avoid this kind of functions in this paper. 
\end{remark}

We now derive an explicit error bound for approximating a function in the regular family by a two-layer neural network with sigmoid functions as the activation functions.  
 \begin{theorem}\label{thm:kai}
     Let $\sigma$ be the standard sigmoid function 
     \begin{equation}
         \sigma(x) = \frac{1}{1+e^{-x}}
     \end{equation}
     then for any $f$ in the regular family and $\varepsilon>0$, there exists a two layer neural network with hidden units $s=\mathcal{O}\left( \frac{d^2}{\varepsilon^2} \right)$~(the first layer) and $r=\mathcal{O}\left( \frac{d}{\varepsilon} \right)$~(the second layer), and constants $d_i$, $c_{ij}$, $\theta_{ij}$, $\gamma_i$ and vectors $\bw^{ij}\in \RR^d$ for which 
    \begin{equation}\label{equ:main}
  \left| {f({\mathbf{x}}) - \sum\limits_{i = 1}^{r} {{d_i}\sigma \left( {\sum\limits_{j = 1}^{s} {{c_{ij}}\sigma ({{\mathbf{w}}^{ij}} \cdot {\mathbf{x}} - {\theta _{ij}})}  - {\gamma _i})} \right)} } \right| < \varepsilon 
\end{equation}
for all $\mathbf{x}=(x_1, \ldots, x_d)\in [0,1]^d$.
     
 \end{theorem}
\begin{proof}
    See \cite{2018arXiv181208883X}.
\end{proof}  
 The theory for multiple layer~(greater than 2) neural networks is elusive. However, \Cref{thm:1d,thm:kai} shed lights on the potential advantage of deeper neural networks: the hidden neuron sizes grow exponentially with the input dimensionality for one layer; they grow polynomially for two layers. We conjecture that for sufficiently many layers such as $\mathcal{O}(d)$ layers, the numbers can be further reduced to $\mathcal{O}(1)$. For example, the famous ResNet has hundreds of layers but each layer is a small convolution filter. The deep but thin neural network enables us to learn high dimensional functions and bypass the curse of the dimensionalities. The investigation will be left to the future. 

\subsection{Error Analysis for Diffusion Equations}

We now conduct error analysis for the framework. To keep the discussion concrete and simple, we consider the 1D diffusion equation where the conductivity is an unknown continuous function. 
\begin{align}\label{equ:model}
    u_t(x, t) =& f(x) u_{xx}(x, t) \qquad x\in [-1,1], t\in (0, T] \\
    u(\pm1, t) = & \varphi_{\pm1}(t)\qquad t\in (0,T] \\
    u(x,0) = & u_0(x)\qquad x\in [-1,1]
\end{align}

The observed dataset is 2 sequential snapshots at time $t=t_0$
\begin{equation}
    M = \{(U_0 = u(\bx, t) + \bW_0, U_1= u(\bx, t) + \bW_1)| \bx = \{x_1, x_2, \ldots, x_n\} \}
\end{equation}
and $\bW_0$ and $\bW_1$ are noise. To put it in another way, given the observed dataset $M\in\RR^{2n}$ and the model specification \cref{equ:model}, we want to find the best function $f_\theta$ that can explain the dataset. We discretize \cref{equ:model} using the uniform grids. Let $-1<x_1<x_2<\ldots<x_n=1$ be the spacial grids with $x_{i+1}-x_i=h$. 

Before we give the convergence rate of the neural network approach, we first investigate several error sources under the framework.

\paragraph{Observation error} $\bW_0$ and $\bW_1$ in this example represent \textit{observation error}, which may come from the measurement error and/or errors in model specification. The former is easy to understand. For example, since there is inherently unpredictable fluctuations in the readings of a measurement apparatus or the experimenter's interpretation of the instrumental reading, the \textit{random error} is almost inevitable in the observations. The errors in the model specification are due to the unrealistic assumption in our models. For instance, we have assumed that the underlying stochastic process is the Wiener process. However, it is well-known that for many physical phenomena such as hydrodynamics, quantum mechanics and so on may follow sub-diffusion/super-diffusion in some cases~\cite{chen2010anomalous}. In other words, the observations suffer from \textit{system errors}.

Yet, in this article, we assume that the dataset error is the random error, i.e., the model specification \cref{equ:model} is correct but the observations may have i.i.d. noise $\bW_1$, $\bW_2$. To quantify the error, we make the following assumption
\begin{assumption}[Observation Error]\label{assum:obs}
There exists a positive constant $\varepsilon_d>0$ such that
\begin{equation}
    \|\bW_0\|_\infty \leq \varepsilon_o \qquad \|\bW_1\|_\infty \leq \varepsilon_o
\end{equation}     
\end{assumption}

\paragraph{Consistency Error} Another important source of error comes from the numerical discretization of \cref{equ:model}. We adopt the standard Crank Nicolson scheme for this problem, which has second-order accuracy both in time and space, which reads
\begin{equation}\label{equ:eq}
    \frac{{v_i^{N + 1} - v_i^N}}{{\Delta t}} = {f_i}\frac{{v_{i + 1}^N + v_{i - 1}^N - 2v_i^N}}{{{2h^2}}} + {f_i}\frac{{v_{i + 1}^{N + 1} + v_{i - 1}^{N + 1} - 2v_i^{N + 1}}}{{{2h^2}}}\qquad i=2,3,\ldots,n-1
\end{equation}
together with boundary conditions
\begin{equation}\label{equ:bd}
    v_n^{N+1} = \varphi(1)\qquad v_1^{N+1} = \varphi(-1)
\end{equation}
Note in \cref{equ:eq} we have used one-step discretization. This is because we have adopted a semi-implicit scheme and it is still stable for reasonable large $\Delta t$. In other cases such as an explicit scheme is used or high precision in temporal direction is required, we may consider multiple steps and formulate an accrued loss function. 

\Cref{equ:eq,equ:bd} leads to a tridiagonal linear system. The method is known to have a local error with an upper bound $\varepsilon_c \leq C_1(\Delta t^2 + h^2) = \mathcal{O}(\Delta t^2 + h^2)$  for some $C_1>0$.
\[\left| {\frac{{u_i^{N + 1} - u_i^N}}{{\Delta t}} - {f_i}\frac{{u_{i + 1}^N + u_{i - 1}^N - 2u_i^N}}{{{2h^2}}} - {f_i}\frac{{u_{i + 1}^{N + 1} + u_{i - 1}^{N + 1} - 2u_i^{N + 1}}}{{{2h^2}}}} \right| \leqslant {\varepsilon _c}\qquad i=2,3,\ldots, n-1\]

The loss function \cref{equ:loss} can be then derived easily
\begin{equation}\label{equ:diffloss}
  \mathrm{loss}(\theta) = \sum\limits_{i = 2}^{n - 1} {{{\left| {\frac{{{U_{1i}} - {U_{0i}}}}{{\Delta t}} - {f_\theta }({x_i})\frac{{{U_{1,i + 1}} + {U_{1,i - 1}} - 2{U_{1i}}}}{{{2h^2}}} - {f_\theta }({x_i})\frac{{{U_{0,i + 1}} + {U_{0,i - 1}} - 2{U_{0,i}}}}{{{2h^2}}}} \right|}^2}} 
\end{equation}

Note we have not included the boundary term since it does not include information about $\theta$. In addition, even if $\mathrm{loss}(\theta)$ is minimized to zero, it does not guarantee that we will obtain a good approximation to $f(x)$. Instead, in this case, the system begins to fit the error, which is known as \textit{overfitting}. There are many techniques to alleviate overfitting; for example, we can collect more observations, use regularizers, terminate the optimization early, etc. We will not dive deep into this topic in the paper since for all the numerical examples we have not observed such issues. 

\paragraph{Optimization error} Errors of this type may arise because the neural network architecture has been poorly designed or the inverse problems itself have multiple solutions~(and therefore the loss function is multimodal). One particular problem the community has recognized is the local minimums for the non-convex loss function. However, the community seems very positive about escaping the local minimums using optimization techniques such as stochastic gradient descent~\cite{kleinberg2018alternative} or specific neural network architectures~\cite{li2018visualizing}. In addition, local minimums do not necessarily mean bad results. As we have discussed, driving the loss function to zero may cause overfitting. A local minimum may give a reasonable estimation of $f(x)$ in many cases. In this paper, we will show that if the other sources~(observation error and consistency error) of error is small, a small optimization error will yield a good estimation. To avoid sophisticated discussion on the optimization phase, we make the following simple assumption
\begin{assumption}[Optimization Error]\label{assum:opt}
There exists a positive constant $\varepsilon_{opt}>0$ such that
\[\left| {\frac{{{U_{1i}} - {U_{0i}}}}{{\Delta t}} - {f_\theta }({x_i})\frac{{{U_{1,i + 1}} + {U_{1,i - 1}} - 2{U_{1i}}}}{{{2h^2}}} - {f_\theta }({x_i})\frac{{{U_{0,i + 1}} + {U_{0,i - 1}} - 2{U_{0,i}}}}{{{2h^2}}}} \right| \leqslant {\varepsilon _{opt}}\]
\end{assumption}

For the convenience of discussion, we assume that the time corresponds to $U_0$ and $U_1$ are $t_0$ and $t_1$ respectively. Now we are ready to state our main theorem
\begin{theorem}\label{thm:main}
    If the following assumptions are satisfied
    \begin{itemize}
        \item Assumptions \ref{assum:obs}  and \ref{assum:opt} are satisfied.
        \item According to Assumption \ref{assum:bound} and  \Cref{thm:uniform_bounded}, there exists constants $F_0$, $F_2>0$ such that 
        \begin{equation}
            |f_\theta(x)| \leq F_0, \qquad |f''_\theta(x)|\leq F_2 \qquad \forall x\in [-1,1]
        \end{equation}
    \item The exact solution $u\in C^4([-1,1]\times [t_0, t_1])$ and
    \begin{equation}
        \delta_2 = \arg\min_{(x,t)\in [-1,1]\times [t_0, t_1]}|u_{xx}(x,t)|\qquad \delta_4 = \arg\max_{(x,t)\in [-1,1]\times [t_0, t_1]}|u_{xxxx}(x,t)|
    \end{equation}
    here $\delta_2$ is strictly positive, i.e.,
    \begin{equation}\label{equ:d2}
        \delta_2 > 0
    \end{equation}
    
    The exact conductivity function $f\in C^2([-1,1])$ and denote
    \begin{equation}
        F^f_2 = \max_{x\in [-1,1]}|f''(x)| 
    \end{equation}
    \item $h$ is sufficiently small in the sense that 
    \begin{equation}\label{equ:hh}
        h < \sqrt{\frac{6\delta_2    }{\delta_4}}
    \end{equation}
    \end{itemize}
    Then the calibrated $f_\theta(x)$ has error estimate $\forall x\in [-1+h, 1-h]$
    \begin{equation}
        |f_\theta(x) - f(x)| \leq \frac{{2{C_1}}}{{{\delta _2}}}\Delta {t^2} + \left( {\frac{{2{C_1}}}{{{\delta _2}}} + \frac{{{F_2} + F_2^f}}{2}} \right){h^2} + \frac{2}{{{\delta _2}}}{\varepsilon _{opt}} + \frac{4}{{{\delta _2}}}\left( {\frac{1}{{\Delta t}} + \frac{{2{F_0}}}{{{h^2}}}} \right){\varepsilon _o}
    \end{equation}
    In other words, 
    \begin{equation}
        |f_\theta(x) - f(x)| = \mathcal{O}\left(\Delta t^2 + h^2 + \varepsilon _{opt} +\left(\frac{1}{\Delta t} + \frac{1}{h^2} \right) \varepsilon _o\right)
    \end{equation}
\end{theorem}
Note the assumption \cref{equ:d2} is reasonable since if $u_{xx}(x,t)=0$ in some region of $(x,t)$, the observations will convey little information about $f(x)$ since $f(x)u_{xx}(x,t)\equiv 0$. If $u_{xx}(x,t)$ is too small, the inverse problem will be quite ill-conditioned and we may have difficulty regularizing the neural network. Nevertheless, in the numerical example even though the second order derivatives of our exact solution vanish near $x=0$, we are still able to infer $f(0)$ quite accurately. This also demonstrates the robustness of our method. 

\begin{proof}
    The proof is split into two parts. We use the notation $u_{0i}=u(x_i, t_0)$ and $u_{1i}=u(x_i, t_1)$.
    \begin{itemize}
        \item Error bound on $|f_\theta(x_i) - f(x_i)|$.
        
        Let $e_{opt,i}$ be the optimization error
        \begin{equation}\label{equ:p2}
  \frac{{{U_{1i}} - {U_{0i}}}}{{\Delta t}} - {f_\theta }({x_i})\frac{{{U_{1,i + 1}} + {U_{1,i - 1}} - 2{U_{1i}}}}{{{2h^2}}} - {f_\theta }({x_i})\frac{{{U_{0,i + 1}} + {U_{0,i - 1}} - 2{U_{0,i}}}}{{{2h^2}}} = {e_{opt,i}}
\end{equation}
        Note that $U_{0i} = u_{0i} + \bW_{0i}$ and $U_{1i} = u_{1i} + \bW_{1i}$, plug them into \cref{equ:p2} we have
        \begin{equation}\label{equ:p3}
  \begin{aligned}
            & \frac{{{u_{1i}} - {u_{0i}}}}{{\Delta t}} - {f_\theta }({x_i})\frac{{{u_{1,i + 1}} + {u_{1,i - 1}} - 2{u_{1i}}}}{{{2h^2}}} - {f_\theta }({x_i})\frac{{{u_{0,i + 1}} + {u_{0,i - 1}} - 2{u_{0,i}}}}{{{2h^2}}} \\
            = & {e_{opt,i}} - \left( {\frac{{{{\mathbf{W}}_{1i}} - {{\mathbf{W}}_{0i}}}}{{\Delta t}} - {f_\theta }({x_i})\frac{{{{\mathbf{W}}_{1,i + 1}} + {{\mathbf{W}}_{1,i - 1}} - 2{{\mathbf{W}}_{1i}}}}{{{2h^2}}} - {f_\theta }({x_i})\frac{{{{\mathbf{W}}_{0,i + 1}} + {{\mathbf{W}}_{0,i - 1}} - 2{{\mathbf{W}}_{0,i}}}}{{{2h^2}}}} \right) 
        \end{aligned}
\end{equation}

        Let $e_{c,i}$ be the consistency error at $x_i$ we have
        \begin{equation}\label{equ:p4}
  \frac{{{u_{1i}} - {u_{0i}}}}{{\Delta t}} - f({x_i})\frac{{{u_{1,i + 1}} + {u_{1,i - 1}} - 2{u_{1i}}}}{{{2h^2}}} - f({x_i})\frac{{{u_{0,i + 1}} + {u_{0,i - 1}} - 2{u_{0,i}}}}{{{2h^2}}} = {e_{c,i}}
\end{equation}

Subtracting \cref{equ:p4} from \cref{equ:p3} we have 
\begin{equation}\label{equ:p8}
  \begin{aligned}
    & \left( {f({x_i}) - {f_\theta }({x_i})} \right)\left[ {\frac{{{u_{1,i + 1}} + {u_{1,i - 1}} - 2{u_{1i}}}}{{2{h^2}}} + \frac{{{u_{0,i + 1}} + {u_{0,i - 1}} - 2{u_{0,i}}}}{{2{h^2}}}} \right] \\
    =&{e_{opt,i}} - {e_{c,i}} \\
    & - \Bigg( {\frac{{{{\mathbf{W}}_{1i}} - {{\mathbf{W}}_{0i}}}}{{\Delta t}} - {f_\theta }({x_i})\frac{{{{\mathbf{W}}_{1,i + 1}} + {{\mathbf{W}}_{1,i - 1}}  - 2{{\mathbf{W}}_{1i}}}}{{2{h^2}}}  - {f_\theta }({x_i})\frac{{{{\mathbf{W}}_{0,i + 1}} + {{\mathbf{W}}_{0,i - 1}} - 2{{\mathbf{W}}_{0,i}}}}{{2{h^2}}}} \Bigg)
\end{aligned}
\end{equation}

Since $\delta_2 > 0$ we must have $u_{xx}\geq \delta_2$ or $u_{xx}\leq \delta_2$ for $(x,t)\in [-1,1]\times [t_0,t_1]$. We might as well assume $u_{xx}\geq \delta_2>0$. 

Note that
\[\frac{{{u_{1,i + 1}} + {u_{1,i - 1}} - 2{u_{1i}}}}{{{h^2}}} = {u_{xx}}({x_i},{t_1}) + \frac{{{h^2}{u_{xxxx}}({x_i},{t_1})}}{{12}}\]
we have
\[\frac{{{u_{1,i + 1}} + {u_{1,i - 1}} - 2{u_{1i}}}}{{{h^2}}} \geqslant {\delta _2} - \frac{{{h^2}{\delta _4}}}{{12}}\]
invoking \cref{equ:hh} we obtain
\[\frac{{{u_{1,i + 1}} + {u_{1,i - 1}} - 2{u_{1i}}}}{{{h^2}}} \geqslant \frac{1}{2}{\delta _2}\]

Likewise
\[\frac{{{u_{0,i + 1}} + {u_{0,i - 1}} - 2{u_{0i}}}}{{{h^2}}} \geqslant \frac{1}{2}{\delta _2}\]
therefore, we have 
\begin{equation}\label{equ:p7}
\frac{{{u_{1,i + 1}} + {u_{1,i - 1}} - 2{u_{1i}}}}{{2{h^2}}} + \frac{{{u_{0,i + 1}} + {u_{0,i - 1}} - 2{u_{0,i}}}}{{2{h^2}}} \geqslant \frac{1}{2}\frac{1}{2}{\delta _2} + \frac{1}{2}\frac{1}{2}{\delta _2} = \frac{1}{2}{\delta _2}
\end{equation}

On the other hand, invoking 
\begin{equation}\label{equ:p6}
    |e_{opt,i}| \leq \varepsilon_{opt}\qquad |e_{c,i}|\leq \varepsilon_c\qquad |W_{ki} |\leq \varepsilon_o\qquad k=0,1
\end{equation}
Combining \cref{equ:p6,equ:p7,equ:p8} we have
\begin{equation}\label{equ:err1}
    \left| {f({x_i}) - {f_\theta }({x_i})} \right| \leqslant \frac{2}{{{\delta _2}}}\left[ {{\varepsilon _{opt}} + {\varepsilon _c} + 2\left( {\frac{1}{{\Delta t}} + \frac{{2{F_0}}}{{{h^2}}}} \right){\varepsilon _o}} \right]
\end{equation}

\item Error bound on $|f_\theta(x) - f(x)|$. 

Define the error function $F:[-1,1]\rightarrow \RR$ as $F(x) = f_\theta(x) - f(x)$, we have already obtained the error bound for $F(x_i)$ in \cref{equ:err1}. For any other point $x\in [-1+h, 1-h]$, we can always find the interval that contains $x$: $x\in [x_i,x_{i+1}]$. Using Taylor's theorem we have
\begin{align}
    F({x_i}) =& F(x) + F'(x)({x_i} - x) + \int_x^{{x_i}} {({x_i} - y)F''(y)dy} \\
    F({x_{i + 1}}) =& F(x) + F'(x)({x_{i + 1}} - x) + \int_x^{{x_{i + 1}}} {({x_{i + 1}} - y)F''(y)dy} 
\end{align}
Let $\alpha  = \frac{{{x_{i + 1}} - x}}{h}$, then we have
\begin{align}
    & \alpha F({x_i}) + (1 - \alpha )F({x_{i + 1}})  \\
=    & F(x) + \alpha \int_x^{{x_i}} {({x_i} - y)F''(y)dy}  + (1 - \alpha )\int_x^{{x_{i + 1}}} {({x_{i + 1}} - y)F''(y)dy} 
\end{align}
We obtain
\[|F(x)| \leqslant \alpha |F({x_i})| + (1 - \alpha )|F({x_{i + 1}})| + \frac{1}{2}|F''(x){|_\infty }{h^2}\]

Note
\[|F''(x)| = |f'{'_\theta }(x)| + |f''(x)| \leqslant {F_2} + F_2^f\]
we have
\begin{align}
    |F(x)| \leqslant & \frac{2}{{{\delta _2}}}\left[ {{\varepsilon _{opt}} + {\varepsilon _c} + 2\left( {\frac{1}{{\Delta t}} + \frac{{2{F_0}}}{{{h^2}}}} \right){\varepsilon _o}} \right] + \frac{{{h^2}}}{2}({F_2} + F_2^f)\\
     \leqslant & \frac{{2{C_1}}}{{{\delta _2}}}\Delta {t^2} + \left( {\frac{{2{C_1}}}{{{\delta _2}}} + \frac{{{F_2} + F_2^f}}{2}} \right){h^2} + \frac{2}{{{\delta _2}}}{\varepsilon _{opt}} + \frac{4}{{{\delta _2}}}\left( {\frac{1}{{\Delta t}} + \frac{{2{F_0}}}{{{h^2}}}} \right){\varepsilon _o}
\end{align}
which finishes the proof. 

    \end{itemize}
\end{proof}

\begin{remark}
    We can see that discrepancy between the estimation $f_\theta(x)$ and the exact function $f(x)$ is bounded by a linear combination of three error sources $\varepsilon_o$, $\varepsilon_d$ and $\varepsilon_{opt}$. $\varepsilon_o$ is usually small if precise measurement can be conducted. When $\Delta t$ or $h$ is large, the consistency error dominates, and we expect to draw a ``convergence plot'' with respect to $\Delta t$ and $h$. If $\Delta t$ and $h$ decrease to a certain level, the optimization error dominates and we will see the error saturates. It is either because the optimizer got stuck at some local minimum or the theoretical loss minimum is reached~(which is called \textit{Bayes error} in the machine learning community -- a statistical lower bound on the error achievable for a given classification problem and associated choice of features~\cite{tumer2003bayes}). In the first case, it may be possible that we can obtain different loss values for different runs~(with a randomized initial guess). In general, the optimization error will be driven to Bayes error with the improvement of optimization techniques. 
\end{remark}

\subsection{Sensitivity Analysis}\label{sect:sens}

The parametrized function $f_\theta(x)$ offers us a way to do sensitivity analysis. When we are most concerned about a physical quantity, for example, the maximum value of the conductivity at a particular location $x = x^*$ in our case, we want to know how reliable the prediction is under our framework. It is similar to uncertainty quantification, but there are two distinct differences: 
\begin{itemize}
    \item The estimation is anchored to a particular physical quantity instead of considering the solution as a whole. 
    \item The estimation will depend on all three sources of errors as well as the neural network structure. In some sense, the estimation will also provide an approach to accessing the quality of the framework in general.  
\end{itemize}

Assume that the physical quantity is expressed as
\begin{equation}
    q(\theta) = Q(f_\theta) 
\end{equation}
where $Q$ is functional. For example
\begin{equation}\label{equ:qf}
    Q(f_\theta) = \max_{x\in [-1,1]}f_\theta(x)
\end{equation}
It is then possible to see the sensitivity of the quantity $q(\theta)$ with respect to the parameter $\theta$ by taking the gradient 
\begin{equation}
    \nabla_\theta q(\theta) = \nabla_\theta Q(f_\theta)
\end{equation}
the latter can usually be done by first discretization $Q(f_\theta)$ and then perform the automatic differentiation. For example, for \cref{equ:qf} we have
\begin{align}
    \nabla_\theta q(\theta) \approx & \nabla_\theta \max_{i=1, 2,3,\ldots, n}\{f_\theta(x_i) \} \\
    = & \nabla_\theta f_\theta(x_{i^*})\qquad i^* = \arg\max_{i=1, 2,3,\ldots, n}\{f_\theta(x_i) \}
\end{align}

$\nabla_\theta f_\theta(x_{i^*})$ can be easily computed with automatic differentiation. After obtaining $\nabla_\theta q(\theta)$, we can give a sensitive region of the $f_\theta(x)$ by
\begin{equation}\label{equ:sensitivity}
    S_{\theta, q, \delta} = \{ f_{\theta'}(x) | \theta' = \theta + \alpha  \nabla_\theta q(\theta), |\alpha|\leq \delta \}
\end{equation}
where $\delta>0$ is a tunable parameter that quantize the sensitivity. 

The intuition is that if we anchor the physical quantity $q(\theta)$, it is most sensitive to the change of hidden parameters $\theta$ in the positive/negative direction $\nabla_\theta q(\theta)$. Therefore $S_{\theta, q, \delta}$ will be a stripe that contains $f_\theta(x)$~(when $\alpha=0$). If the stripe deviates from the $f_\theta$ too much or changes dramatically for a reasonable $\alpha$, the framework is then questionable, and we should investigate components such as the neural network structure. The sensitivity analysis provides us with a way to evaluate the quality of the framework. 

\begin{figure}[htbp]
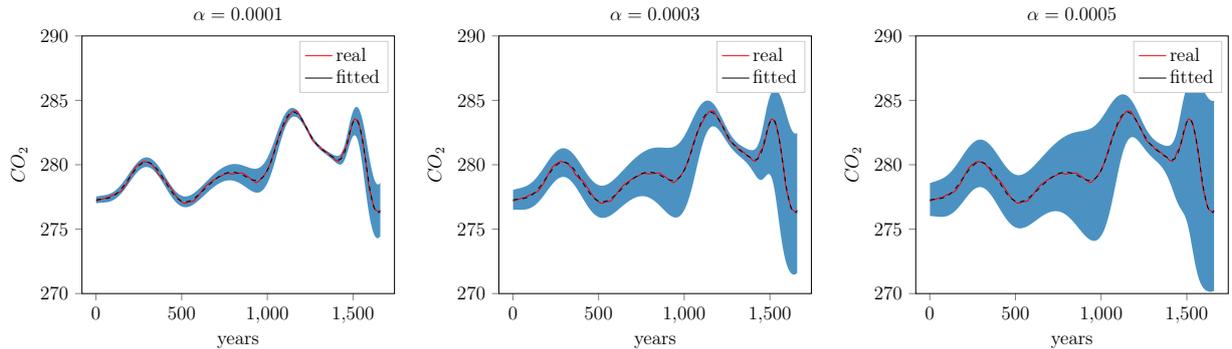

\centering
\scalebox{0.6}{\input{figures/co21.tex}}~
\scalebox{0.6}{\input{figures/co22.tex}}~
\scalebox{0.6}{\input{figures/co23.tex}}
\caption{The width of the stripe describes the impact of a small change in the neural network parameters if the maximum value of the fitted value is the targeted physical quantity.}
\label{fig:CO2}
\end{figure}

As an example, we consider fitting the global yearly mean data~\footnote{\url{ftp://data.iac.ethz.ch/CMIP6/input4MIPs/UoM/GHGConc/CMIP/yr/atmos/UoM-CMIP-1-1-0/GHGConc/gr3-GMNHSH/v20160701/mole_fraction_of_carbon_dioxide_in_air_input4MIPs_GHGConcentrations_CMIP_UoM-CMIP-1-1-0_gr3-GMNHSH_0000-2014.csv}} for some quantity related to $\mathrm{CO}_2$. We are most concerned about the maximum of the data and want to see how sensitive it is with respect to the neural network parameters. \Cref{fig:CO2} shows the sensitivity region computed using the method described above. Here $\alpha$ is the largest step size in \cref{equ:sensitivity}. The width of the stripe describes the impact of a small change in the neural network parameters if the maximum value of the fitted value is the targeted physical quantity. 

The study of the sensitivity analysis is at an early stage. Many problems need to be addressed. For example, how can we pick the reasonable $\alpha$? Also, the relationship between the sensitivity analysis and uncertainty quantification of the model is not clear and requires more theoretical investigation. That will be left to the future work. Nevertheless, the numerical examples demonstrate that it captures the sensitivity of the physical quantity depending on the neural network model, and can potentially serve as a way for the diagnosis of the framework. 

\section{Numerical Examples}

We now apply the framework to three differential equations from various applications and demonstrate its effectiveness. The first equation is a diffusion equation from the elliptic equation family. We try to infer the unknown conductivity function. The second is a wave equation from the hyperbolic equation family. The unknown is the velocity fields. The problem is a significant and challenging one in seismic imaging. Instead of using the popular adjoint method used by the geophysics community, we provide another potential way to infer the velocity fields. The last problem is a variable coefficient nonlinear Burger's equation. The Burger's equation appeared as one of the simplest instances of a nonlinear system out of equilibrium in fluid dynamics~\cite{bec2007burgers}. The study may offer new approaches to the inverse problems in turbulence modeling. 

For all the numerical experiment below, we have the following common settings
\begin{itemize}
    \item For the neural network model, we let $n_l=3$ and the number of hidden neurons per layer be $n_h = 20$. We have tested extensively for different $n_l=3,4,5,6,7,8,9,10$ and $n_h=20,40,60,80$ but find basically no difference except for large $n_l$ and $n_h$, the optimizer has difficulty converging within the predetermined number of iterations. 
    \item We use \texttt{L-BFGS} optimizer, which enjoys convergence property similar to second order method. The iterator will stop if either of the following is satisfied 
    \begin{align}
        |f(\bx_{k+1})-f(\bx_k)|\leq & \varepsilon_1 |f(\bx_k)|\\
        |\nabla f(\bx_{k+1})| \leq & \varepsilon_2
    \end{align}
    where $\bx_{k}$ is the variable value at iteration $k$. We set $\varepsilon_1=10^{-12}$, $\varepsilon_2=10^{-12}$. In addition, if not mentioned, the number of maximum iteration is set as 5000. 
\end{itemize}

\subsection{Diffusion Equation}

We let the true model be
\begin{equation}\label{equ:ex1}
    u_t(x,t) = c(x) u_{xx}(x,t) + f(x)\qquad (x,t)\in [-1,1]\times [0,T]
\end{equation}
where $u(x) = e^{-\pi^2 t}\sin(\pi x)$, $c(x) = 1+e^{-(x-0.5)^2}$ and
\begin{equation*}
f(x,t) = \pi^{2} \left(1.0 + e^{- \left(x - 0.5\right)^{2}}\right) e^{- \pi^{2} t} \sin{\left (\pi x \right )} - \pi^{2} e^{- \pi^{2} t} \sin{\left (\pi x \right )}
\end{equation*}

\begin{figure}[htbp]
\centering
\scalebox{1.0}{\input{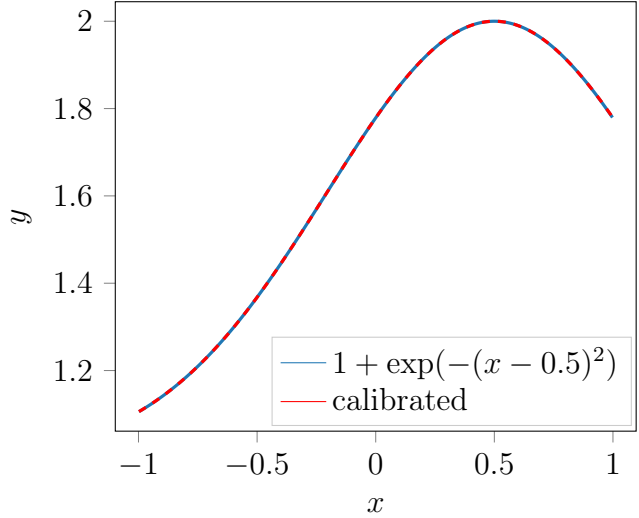}}
\caption{The snapshots are taken at $t=0.1$ and $t=t+\Delta t$, with $\Delta t=0.001$ and $h=0.002$. We can see that the calibrated value of $c(x)$ matches perfectly with the exact value.}
\label{fig:pydiff0}
\end{figure}

\Cref{fig:pydiff0} shows the case where the snapshot is taken at $t=0.1$, with $\Delta t=0.001$ and $h=0.002$. We can see that the calibrated value of $c(x)$ matches perfectly with the exact value. In this case, we do not introduce any noise to the solution, i.e., $\varepsilon_o\equiv 0$.

\begin{figure}[htbp]
\centering
\scalebox{1.0}{
\begin{tikzpicture}

\definecolor{color0}{rgb}{0.12156862745098,0.466666666666667,0.705882352941177}
\definecolor{color1}{rgb}{1,0.498039215686275,0.0549019607843137}
\definecolor{color2}{rgb}{0.172549019607843,0.627450980392157,0.172549019607843}
\definecolor{color3}{rgb}{0.83921568627451,0.152941176470588,0.156862745098039}
\definecolor{color4}{rgb}{0.580392156862745,0.403921568627451,0.741176470588235}
\definecolor{color5}{rgb}{0.549019607843137,0.337254901960784,0.294117647058824}
\definecolor{color6}{rgb}{0.890196078431372,0.466666666666667,0.76078431372549}

\begin{axis}[
legend cell align={left},
legend entries={{$n=10$},{$n=20$},{$n=40$},{$n=80$},{$n=160$},{$n=320$},{$n=640$},{$\mathcal{O}(\Delta t^{2})$}},
legend style={at={(1.03,0.2)}, anchor=south west, draw=white!80.0!black},
tick align=outside,
tick pos=left,
x dir=reverse,
title={$n_h=20, n_l=3$},
x grid style={white!69.01960784313725!black},
xlabel={$\Delta t$},
xmin=0.000794328234724282, xmax=0.125892541179417,
xmode=log,
y grid style={white!69.01960784313725!black},
ylabel={Error},
ymin=7.18647984384873e-05, ymax=0.103076473404677,
ymode=log
]
\addlegendimage{no markers, color0}
\addlegendimage{no markers, color1}
\addlegendimage{no markers, color2}
\addlegendimage{no markers, color3}
\addlegendimage{no markers, color4}
\addlegendimage{no markers, color5}
\addlegendimage{no markers, color6}
\addlegendimage{no markers, white!49.80392156862745!black}
\addplot [semithick, color0, mark=*, mark size=3, mark options={solid}]
table [row sep=\\]{%
0.1	0.03670271109645 \\
0.05	0.0463068007282523 \\
0.01	0.0659446885226462 \\
0.005	0.0666068145132439 \\
0.001	0.0667870863580617 \\
};
\addplot [semithick, color1, mark=*, mark size=3, mark options={solid}]
table [row sep=\\]{%
0.1	0.065142785130901 \\
0.05	0.0105159443170648 \\
0.01	0.0157221559908587 \\
0.005	0.0163280000953794 \\
0.001	0.0165035092739889 \\
};
\addplot [semithick, color2, mark=*, mark size=3, mark options={solid}]
table [row sep=\\]{%
0.1	0.0717678975897289 \\
0.05	0.0174981370341978 \\
0.01	0.00330849077036977 \\
0.005	0.00392305155371853 \\
0.001	0.0041095006545282 \\
};
\addplot [semithick, color3, mark=*, mark size=3, mark options={solid}]
table [row sep=\\]{%
0.1	0.0734245641759934 \\
0.05	0.0192325759674647 \\
0.01	0.000226326295597268 \\
0.005	0.000859959404116761 \\
0.001	0.00103793535203911 \\
};
\addplot [semithick, color4, mark=*, mark size=3, mark options={solid}]
table [row sep=\\]{%
0.1	0.0738557918243619 \\
0.05	0.0197139240112187 \\
0.01	0.000752069511007925 \\
0.005	0.000599429383797645 \\
0.001	0.00026587214879803 \\
};
\addplot [semithick, color5, mark=*, mark size=3, mark options={solid}]
table [row sep=\\]{%
0.1	0.0739495862355244 \\
0.05	0.0199453598639545 \\
0.01	0.000874095004519182 \\
0.005	0.000424145727620839 \\
0.001	0.000511169176337978 \\
};
\addplot [semithick, color6, mark=*, mark size=3, mark options={solid}]
table [row sep=\\]{%
0.1	0.074075699849772 \\
0.05	0.0199938805545525 \\
0.01	0.00105405407183357 \\
0.005	0.000916810946652413 \\
0.001	0.000723232136657925 \\
};
\addplot [semithick, white!49.80392156862745!black, dashed]
table [row sep=\\]{%
0.1	0.01 \\
0.05	0.0025 \\
0.01	0.0001 \\
};
\end{axis}

\end{tikzpicture}}
\scalebox{1.0}{
\begin{tikzpicture}

\definecolor{color0}{rgb}{0.12156862745098,0.466666666666667,0.705882352941177}
\definecolor{color1}{rgb}{1,0.498039215686275,0.0549019607843137}
\definecolor{color2}{rgb}{0.172549019607843,0.627450980392157,0.172549019607843}
\definecolor{color3}{rgb}{0.83921568627451,0.152941176470588,0.156862745098039}
\definecolor{color4}{rgb}{0.580392156862745,0.403921568627451,0.741176470588235}
\definecolor{color5}{rgb}{0.549019607843137,0.337254901960784,0.294117647058824}

\begin{axis}[
legend cell align={left},
legend entries={{$\Delta t=0.1$},{$\Delta t=0.05$},{$\Delta t=0.01$},{$\Delta t=0.005$},{$\Delta t=0.001$},{$\mathcal{O}(h^{2})$}},
legend style={at={(1.03,0.2)}, anchor=south west, draw=white!80.0!black},
tick align=outside,
tick pos=left,
title={$n_h=20, n_l=3$},
x grid style={white!69.01960784313725!black},
x dir=reverse,
xlabel={$h$},
xmin=0.00253828873861328, xmax=0.246228882668983,
xmode=log,
y grid style={white!69.01960784313725!black},
ylabel={Error},
ymin=0.000169429089799733, ymax=0.0989515953877269,
ymode=log
]
\addlegendimage{no markers, color0}
\addlegendimage{no markers, color1}
\addlegendimage{no markers, color2}
\addlegendimage{no markers, color3}
\addlegendimage{no markers, color4}
\addlegendimage{no markers, color5}
\addplot [semithick, color0, mark=*, mark size=3, mark options={solid}]
table [row sep=\\]{%
0.2	0.03670271109645 \\
0.1	0.065142785130901 \\
0.05	0.0717678975897289 \\
0.0249999999999999	0.0734245641759934 \\
0.0125	0.0738557918243619 \\
0.00625000000000009	0.0739495862355244 \\
0.00312500000000004	0.074075699849772 \\
};
\addplot [semithick, color1, mark=*, mark size=3, mark options={solid}]
table [row sep=\\]{%
0.2	0.0463068007282523 \\
0.1	0.0105159443170648 \\
0.05	0.0174981370341978 \\
0.0249999999999999	0.0192325759674647 \\
0.0125	0.0197139240112187 \\
0.00625000000000009	0.0199453598639545 \\
0.00312500000000004	0.0199938805545525 \\
};
\addplot [semithick, color2, mark=*, mark size=3, mark options={solid}]
table [row sep=\\]{%
0.2	0.0659446885226462 \\
0.1	0.0157221559908587 \\
0.05	0.00330849077036977 \\
0.0249999999999999	0.000226326295597268 \\
0.0125	0.000752069511007925 \\
0.00625000000000009	0.000874095004519182 \\
0.00312500000000004	0.00105405407183357 \\
};
\addplot [semithick, color3, mark=*, mark size=3, mark options={solid}]
table [row sep=\\]{%
0.2	0.0666068145132439 \\
0.1	0.0163280000953794 \\
0.05	0.00392305155371853 \\
0.0249999999999999	0.000859959404116761 \\
0.0125	0.000599429383797645 \\
0.00625000000000009	0.000424145727620839 \\
0.00312500000000004	0.000916810946652413 \\
};
\addplot [semithick, color4, mark=*, mark size=3, mark options={solid}]
table [row sep=\\]{%
0.2	0.0667870863580617 \\
0.1	0.0165035092739889 \\
0.05	0.0041095006545282 \\
0.0249999999999999	0.00103793535203911 \\
0.0125	0.00026587214879803 \\
0.00625000000000009	0.000511169176337978 \\
0.00312500000000004	0.000723232136657925 \\
};
\addplot [semithick, color5, dashed]
table [row sep=\\]{%
0.2	0.02 \\
0.1	0.00500000000000001 \\
0.05	0.00125 \\
0.0249999999999999	0.000312499999999998 \\
};
\end{axis}

\end{tikzpicture}}
\caption{We take the snapshots at $t=0.1$ and $t=t+\Delta t$, with various $\Delta t$ and $h$, and $n = \frac{2}{h}$ in the plots. We see a clear pattern exactly predicted in \Cref{thm:main}: a second order convergence with respect to $\Delta t$ and $h$, and at some point, the error ceases to decrease and remains at a constant level. }
\label{fig:pydiff12}
\end{figure}
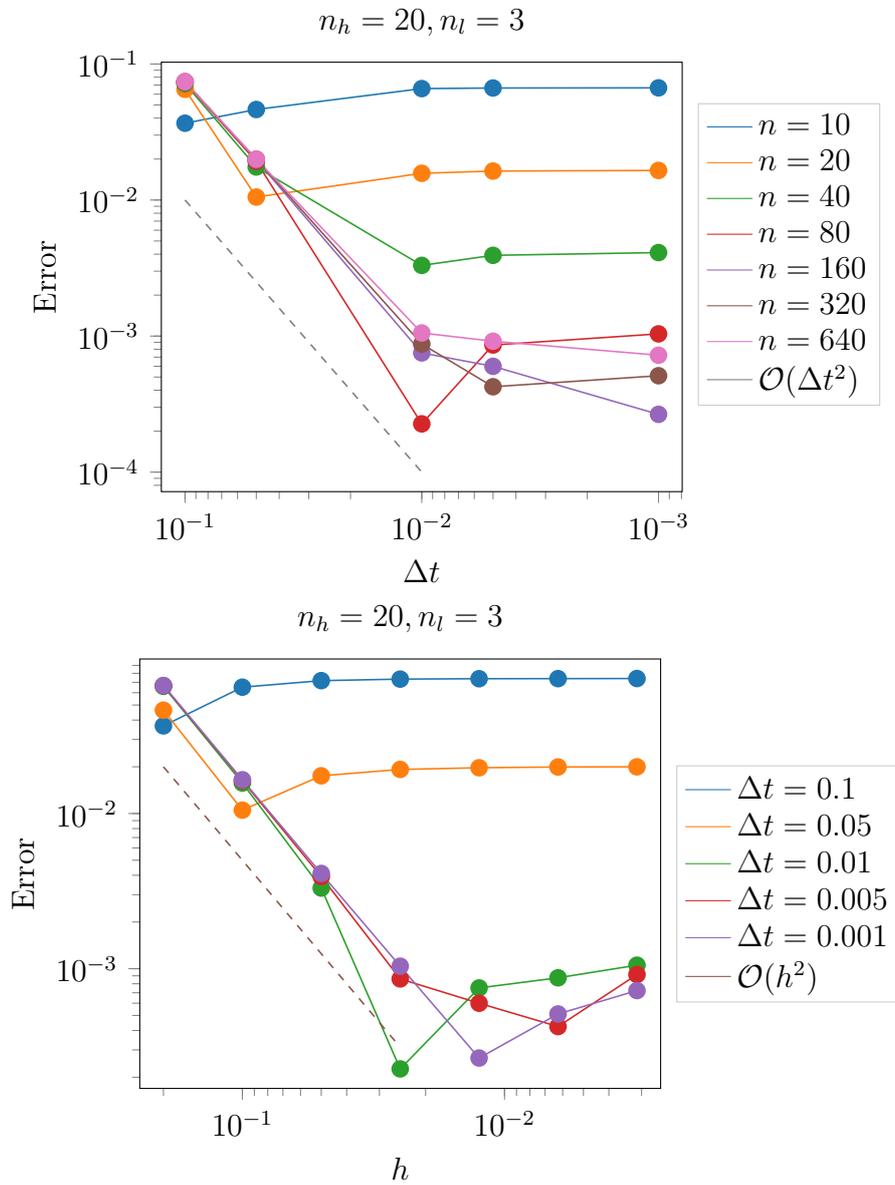

We now conduct a systematic study on the influence of $\Delta t$ and $h$. In \cref{fig:pydiff12} we take the snapshots at $t=0.1$ and $t=t+\Delta t$, with various $\Delta t$ and $h$, and $n = \frac{2}{h}$ in the plots. We see a clear pattern exactly predicted in \Cref{thm:main}: a second-order convergence in $\Delta t$ and $h$, and at some point, the error ceases to decrease and remains at a constant level. The numerical results not only demonstrate our theory but also offer us the direction for further optimization: \textbf{collecting more data is of little benefit after the error saturates; instead, we should focus on optimization or noise reduction.} It is quite a surprising implication since typically people believe more data is beneficial. 

We can also see that when the error saturates, the level of error are not the same for different $n$ or $\Delta t$. This implies that the spacial discretization error in the first plot or the temporal discretization error in the second plot dominates. If $h$ or $\Delta t$ is sufficiently small, the level where the error saturates will be approximately the same, as shown in \cref{fig:pydiff12}.

Finally, we carry out the sensitivity analysis. The settings are the same as that in \cref{fig:pydiff0} except that we add i.i.d. Gaussian noise with mean 0 and standard deviation of $3\times 10^{-7}$ to the observations. The physical quantity we are interested is the conductivity function value at $x=0$, i.e., the prediction for $c(0)$. Following the discussion in \Cref{sect:sens}, we generate the sensitive region of this quantity. We see that the stripe gradually grows and eventually incorporates the exact solution. We also see that the variance near the boundary is smaller than that in the center. This observation is consistent with our intuition that if the interesting physical quantity is the value of $c(x)$ at $0$, the far-away region should be less impacted. Therefore a modification to the neural network parameter will not result in significant change in the far-away value of the calibrated function.

\begin{figure}[htbp]
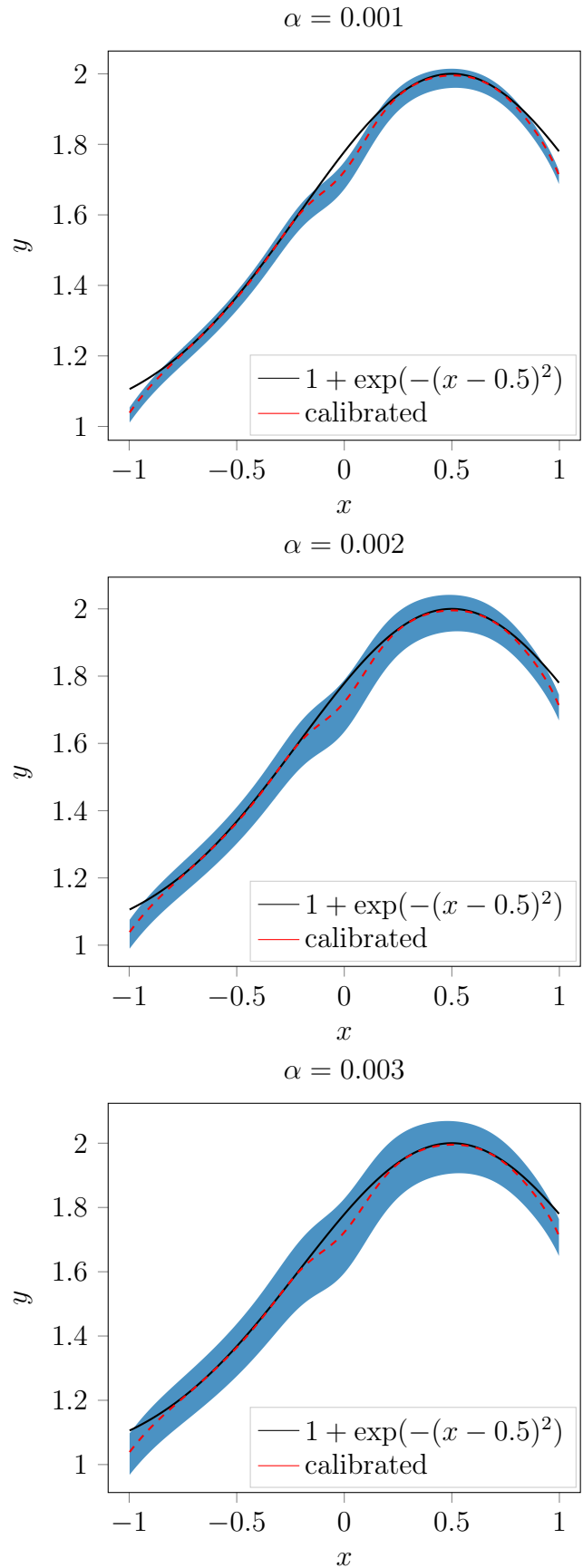

\centering
\scalebox{1.0}{\input{figures/pydiff_uq1.tex}}
\scalebox{1.0}{\input{figures/pydiff_uq2.tex}}
\scalebox{1.0}{\input{figures/pydiff_uq3.tex}}

\caption{Sensitivity analysis for \cref{equ:ex1}. The physical quantity we are interested is the conductivity function value at $x=0$, i.e., the prediction for $c(0)$. We have used three different step size $\alpha=0.001$, $0.002$ and $0.003$. The blue region is the sensitivity region defined by \cref{equ:sensitivity}.}
\label{fig:pydiff_uq}
\end{figure}

\subsection{Wave Equation}

The second example is concerned about wave equations. Consider the model problem
\begin{equation}
  u_{tt}(x, t) = c(x) u_{xx}(x,t)\qquad (x,t)\in [-1,1]\times [0,1]
\end{equation}
where $c(x)$ is the unknown velocity field. We set the initial condition as
\begin{equation}
    u(x,0) = e^{-10x^2}
\end{equation}
The two manufactured velocity fields are
\begin{equation}\label{equ:cc}
    c_1(x) = 1 + \exp(-(x-0.5)^2)\qquad c_2(x) = \begin{cases}
        2 & x\in [-1,-0.5)\\
        1 & x\in [-0.5,-0.1)\\
        0.2 & x\in [-0.1,0.1)\\
        1.0 & x\in [0.2,0.3)\\
        1.5 & x\in (0.3,1.0]
    \end{cases}
\end{equation}

Here we showcase two modifications to the framework and demonstrate its flexibility for practical problems
\begin{itemize}
    \item Suppose we have snapshots at different times and three sequential snapshots for each, i.e., $|\mathcal{T}|>1$. We want to make full use of the observations. Therefore we construct an ensemble loss function as shown in \cref{equ:loss}.
    \item Suppose we have prior knowledge $c(x)\in [0,2]$. To incorporate the prior knowledge, we let the last layer in \cref{equ:nn} be
    \begin{equation}
        \by_{n_l+1}(\by_{n_l}) = \tanh(\bW_{n_l} \by_{n_l} + \bb_{n_l}) + 1
    \end{equation}
\end{itemize}

\begin{figure}[htbp]
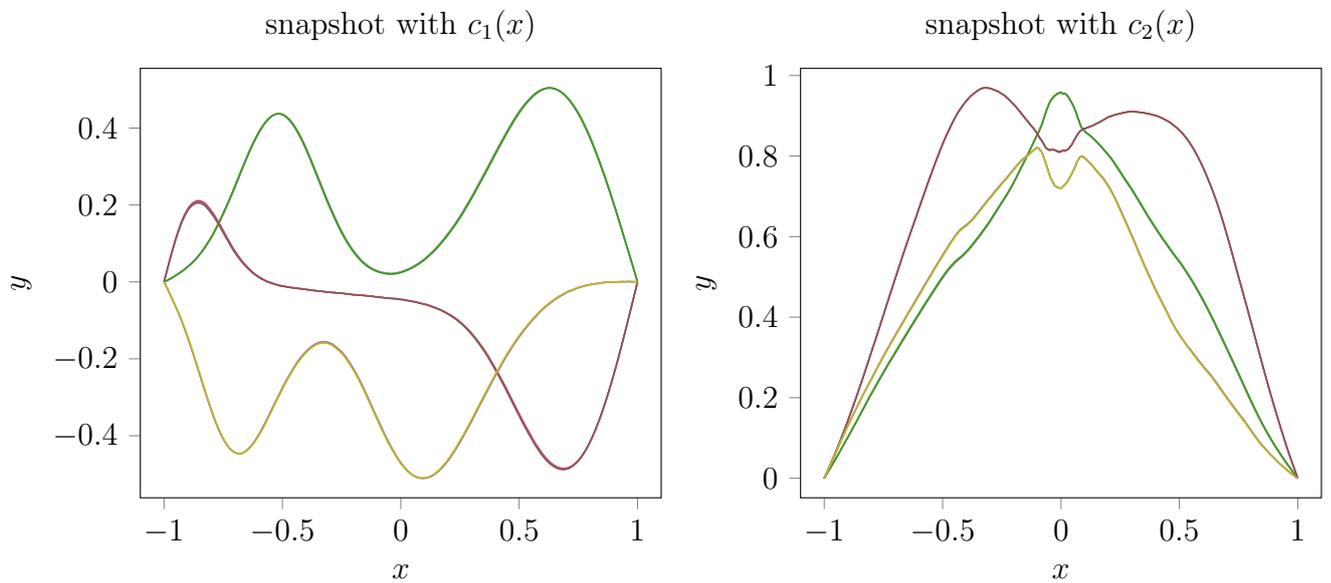

\centering
\scalebox{1.0}{\input{figures/wavetrue_3.tex}}~
\scalebox{1.0}{\input{figures/wavetrue_1.tex}}
\caption{Examples of the snapshots for the velocity fields $c_1(x)$ and $c_2(x)$. Note that since $10\Delta t$ is quite small, the nearby sequential snapshots almost overlap.}
\label{fig:wave}
\end{figure}

Since we may not have analytical solution for the velocity fields \cref{equ:cc}, we first run the simulation with $h = 0.004$, $\Delta t = 0.0001$ using the numerical scheme
\begin{equation}\label{equ:wave}
  \frac{{v_i^{N + 2} - 2v_i^{N + 1} + v_i^N}}{{\Delta {t^2}}} = c_i^2\frac{{v_{i + 1}^{N + 1} - 2v_i^{N + 1} + v_{i - 1}^{N + 1}}}{{\Delta {t^2}}}\qquad i = 2,3,\ldots, n-1
\end{equation}
then we use the values at $t = \frac{k}{3}-20\Delta t$, $\frac{k}{3}-10\Delta t$, $\frac{k}{3}$, $k=1$, $2$, $3$ as the observation. The loss function \cref{equ:loss} is also formulated using \cref{equ:wave} except that the time step is now $10\Delta t$. Note that $\frac{2\cdot 10\Delta t}{h}=0.5<1$, the CFL condition is satisfied. We also add i.i.d. Gaussian noise sampled from $\mathcal{N}(0, (10\Delta t^2)^2)$ to the observations. \Cref{fig:wave} shows examples of the snapshots for the velocity fields $c_1(x)$ and $c_2(x)$. Note that since $10\Delta t$ is quite small, the nearby sequential snapshots almost overlap. However, our algorithm is able to infer the velocity fields from the subtle difference.

\begin{figure}[htbp]
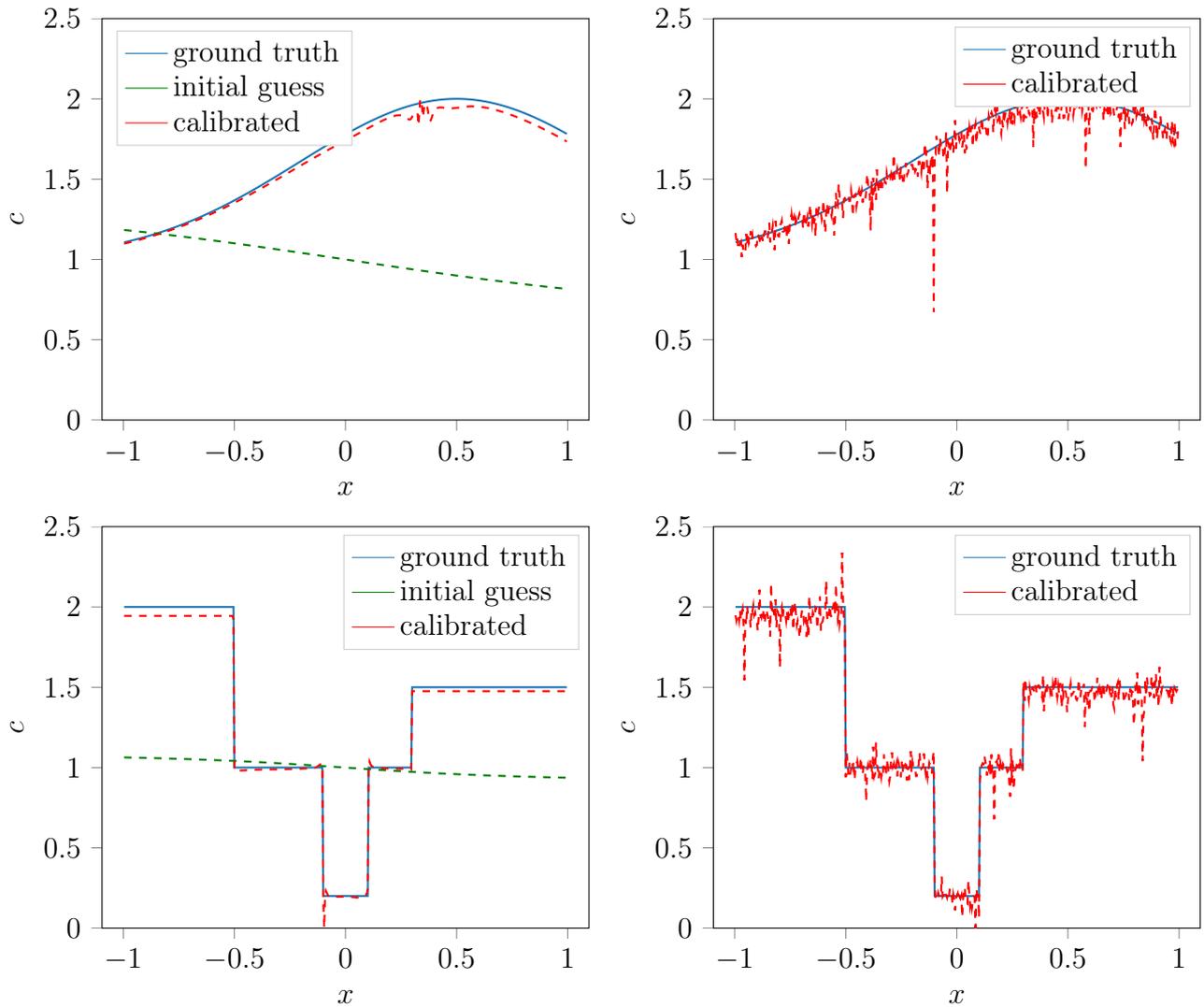

\centering
\scalebox{1.0}{\input{figures/lwavetrue_3.tex}}~
\scalebox{1.0}{\input{figures/lwavefalse_3.tex}}
\scalebox{1.0}{\input{figures/lwavetrue_1.tex}}~
\scalebox{1.0}{\input{figures/lwavefalse_1.tex}}
\caption{The left columns show the calibrating results using the neural networks while the right columns correspond to least square methods. We can see that the least square methods are vulnerable to noise and do not yield a solution as smooth as the neural network methods. Especially for $c_2(x)$, the neural network captures the transition between different velocity levels and yields a constant value at each level. The results are quite remarkable since it is generally challenging to capture those discontinuities in inverse problems. The fit is not perfect though because we have added noise to the dataset.}
\label{fig:lwave}
\end{figure}

We also compare the method with neural networks with that of least square methods. In \cref{fig:lwave}, the left columns show the calibrating results using the neural networks while the right columns correspond to least square methods. For the least square methods, we have used \texttt{Newton-CG} algorithms for minimizing the loss function with discrete velocity field values as unknown variables. We can see that the least square methods are vulnerable to noise and do not yield a solution as smooth as the neural network methods. Especially for $c_2(x)$, the neural network captures the transition between different velocity levels and yields a constant value at each level. The results are quite remarkable since it is generally challenging to capture those discontinuities in inverse problems. The fit is not perfect though because we have added noise to the dataset. 

\subsection{Variable Coefficient Burger's Equation}

Finally, we study a variable-coefficient Burgers equation arising in the modeling of segregation of dry bidisperse granular mixtures~\cite{christov2017numerical}. The numerical PDE is due to \cite{dolgunin1995vn}
\begin{equation}\label{equ:beq}
    \frac{\partial u}{\partial t} + \frac{\partial}{\partial x}(u(1-u)f(x)) - D\frac{\partial^2 u}{\partial x^2} = 0,\qquad x\in [-1,1], t>0
\end{equation}
here $u=u(x,t)$ is the concentration. The function $f(x)$ is called the percolation velocity and $D$ is the constant diffusivity due to inter-particle collisions. In this numerical example, we consider the case where we have partially observed evolution of $u(x,t)$ and want to infer the percolation velocity $f(x)$. To formulate the loss function, we use the following numerical discretization scheme from \cite{christov2017numerical}
\begin{multline}\label{equ:burger}
    \left[ {\frac{{\Delta t}}{{8h}}\left( {2 - v _{j + 1}^{N + 1} - 2v _{j + 1}^N} \right)f({x_{j + 1}}) - \frac{{D\Delta t}}{{2{h^2}}}} \right]v _{j + 1}^{N + 1} + \left[ {1 + \frac{{D\Delta t}}{{{h^2}}}} \right]v _j^{N + 1}\\
     + \left[ { - \frac{{\Delta t}}{{8h}}\left( {2 - v _{j - 1}^{N + 1} - 2v _{j - 1}^N} \right)f({x_{j - 1}}) - \frac{{D\Delta t}}{{2{h^2}}}} \right]v _{j - 1}^{N + 1}\\
      =  - \frac{{\Delta t}}{{8h}}\left[ {v _{j + 1}^N(2 - v _{j + 1}^N)f({x_{j + 1}}) - v _{j - 1}^N(2 - v _{j - 1}^N)f({x_{j - 1}})} \right] + v _j^N\\
       + \frac{{D\Delta t}}{{2{h^2}}}\left( {v _{j + 1}^N - 2v _j^N + v _{j - 1}^N} \right)\qquad j = 2,3,\ldots, n-1
\end{multline}

Note \cref{equ:burger} is a nonlinear equation, and in practice, we can use Newton's methods to obtain $v^{N+1}$ from $v^N$. For the inverse problem, the loss function is constructed by taking the square sum of the difference between the left-hand side and the right-hand side in \cref{equ:burger} over $j = 2$,$3$,$\ldots$, $n-1$.

\begin{figure}[htbp] 
\centering
\includegraphics[width=1.0\textwidth,keepaspectratio]{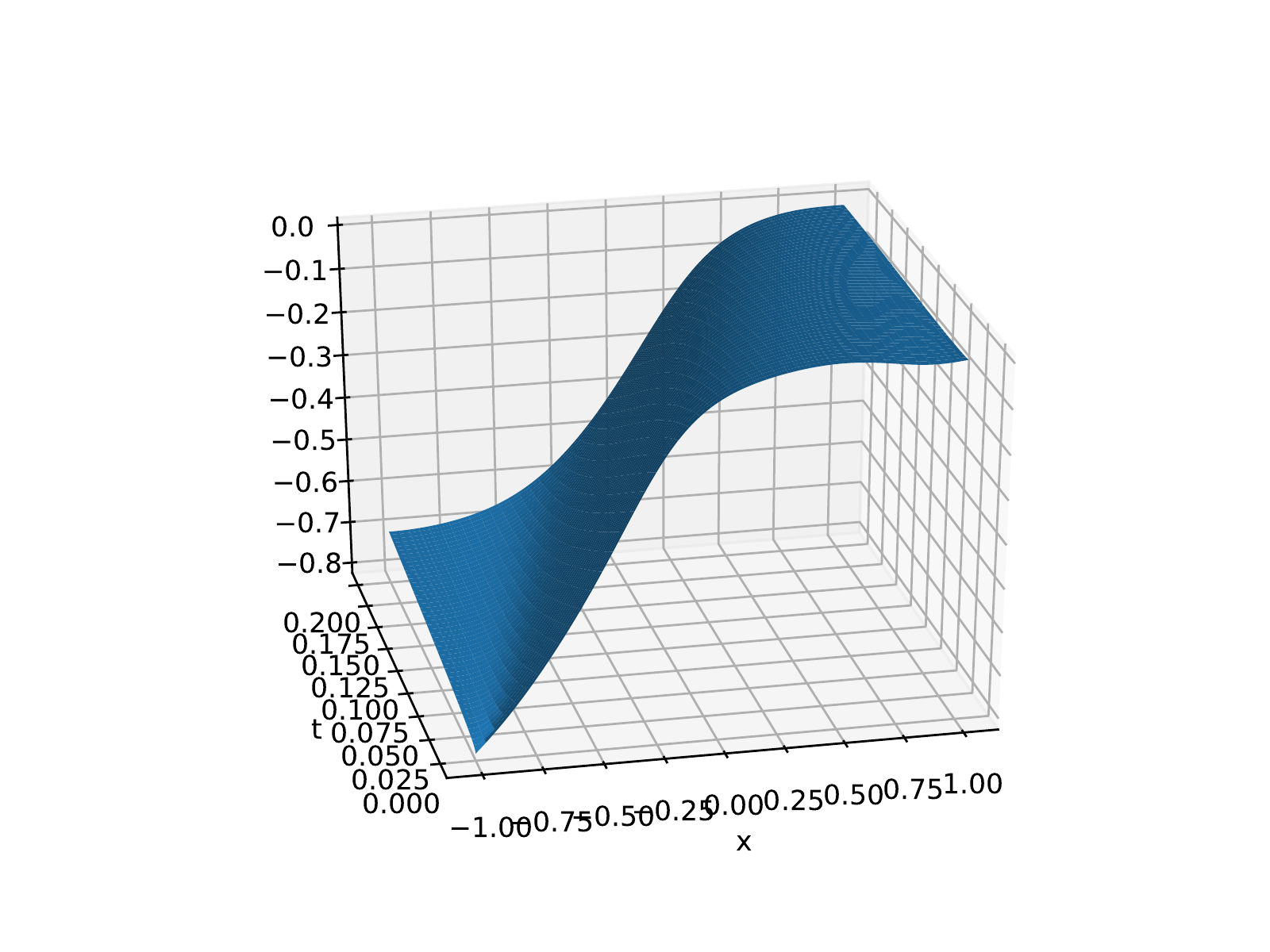}
\caption{Solution profile for \cref{equ:beq} with initial condition \cref{equ:binit} and the percolation function \cref{equ:bfx}.}
\label{fig:bsol}
\end{figure}

For this problem, we use a neural network with four hidden layers and each layer has 40 neurons. We have found that the performance is insensitive to these hyper-parameters as long as they stay in a reasonable range~(neither too deep nor too heavy). The exact solution is simulated using $\Delta t=2\times 10^{-6}$ and $h = 0.004$ while we take 20 snapshots at $t=0.01$,  $0.02$, $0.03$, $\ldots$, $0.20$. Noise is sampled from i.i.d. Gaussian distribution with zero mean and $10^{-5}$ standard deviation. The time step is quite coarse, but we show it does not compromise our algorithm to recover $f(x)$. We fix $D=0.1$. We also apply $L_2$ normalization to the neural networks with penalty $0.01$ to enforce smoothness.  The initial condition is
\begin{equation}\label{equ:binit}
    u(x,0) = \frac{c-f(x)}{2f(x)}\left( -1+\frac{(x-ct)\tanh(c-f(x))}{2D} \right)
\end{equation}
where
\begin{equation}\label{equ:bfx}
    f(x) = -1 + \exp(-(x-0.5)^2)
\end{equation}
The corresponding solution profile is shown in \cref{fig:bsol}.

\begin{figure}[htbp]
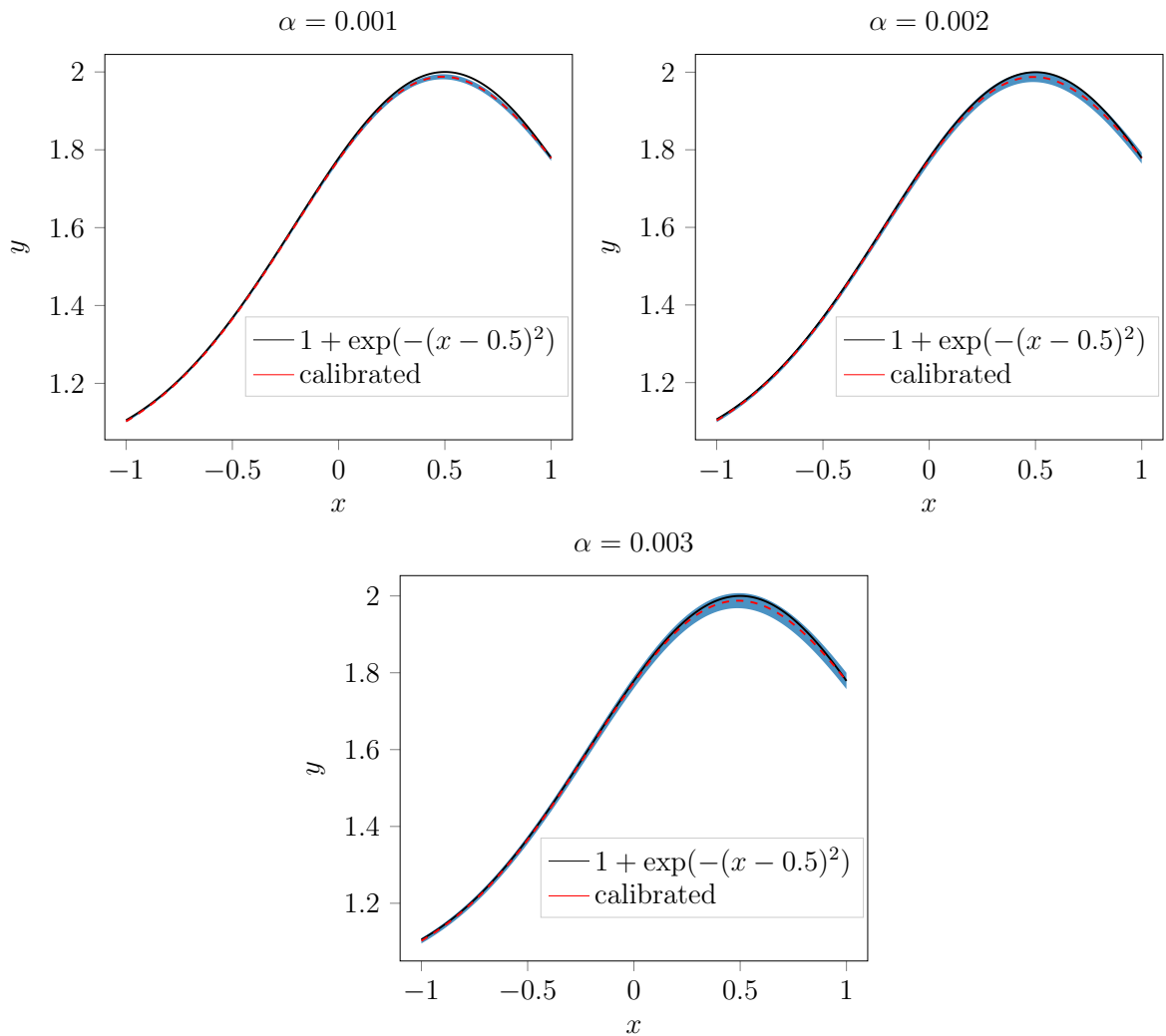

\centering
\scalebox{0.9}{\input{figures/burger_uq1.tex}}
\scalebox{0.9}{\input{figures/burger_uq2.tex}}
\scalebox{0.9}{\input{figures/burger_uq3.tex}}
\caption{The calibrated results and the corresponding sensitivity regions with different step sizes $\alpha$. The stripe gradually grows and finally encloses the exact $f(x)$. We also see that the values far-away from the maxima are less sensitive to the parameters. }
\label{fig:burger}
\end{figure}

For sensitivity analysis, we consider the maximum value of $f(x)$ as the interesting physical quantity, i.e.,
\begin{equation}
    Q(f_\theta) = \max_{x\in [-1,1]}f_\theta(x)
\end{equation}
\Cref{fig:burger} shows the calibrated results and the corresponding sensitivity regions with different step sizes $\alpha$. The stripe gradually grows and finally encloses the exact $f(x)$. We also see that the values far-away from the maxima are less sensitive to the parameters. 

\begin{figure}[htbp]
\centering
\scalebox{1.0}{\input{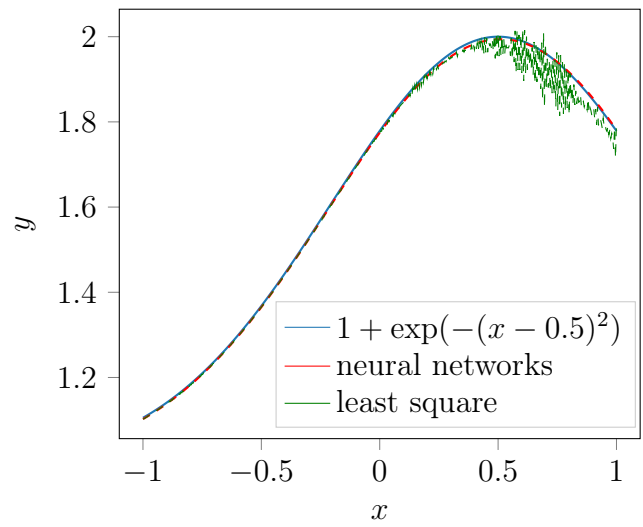}}
\caption{The calibrated results using least square methods with the same settings. We can see that this method is more susceptible to random noise and has stability issues compared to the neural network approach.}
\label{fig:burger0}
\end{figure}

As a comparison, \cref{fig:burger0} shows the calibrated results using least square methods with the same settings. We can see that this method is more susceptible to random noise and has stability issues compared to the neural network approach. 

\section{Conclusion}

We have introduced a new framework for inverse problems in differential equations based on neural networks and automatic differentiation. It leverages the current development of machine learning techniques, especially deep learning. Thanks to the open source frameworks developed by the machine learning community, such as TensorFlow~\cite{tensorflow2015-whitepaper}, PyTorch~\cite{paszke2017automatic}, MXNet~\cite{chen2015mxnet}, and so on, the inverse problem solvers can be easily implemented and incorporated into traditional numerical PDE codes. In future works, we show that the approach can also be easily combined with more sophisticated numerical methods such as finite element methods and finite volume methods, which are widely deployed in software for fluid dynamics, solid mechanics, and more. In this framework, a crucial step is to represent the unknown fields using neural networks, instead of discrete values or a linear combination of kernel basis functions. The procedure offers the universal approximation ability of neural networks, along with regularizing solution smoothness, bypassing the curse of dimensionality and leveraging the high-performance computing developed for deep learning. It works seamlessly with forward-simulation schemes and benefits from the numerical stability and consistency. Also, we are also able to provide convergence analysis and contain the error, using theoretical tools from classical computational mathematics. 

We performed initiatory demonstrations for the effectiveness of the framework on elliptic and hyperbolic problems, as well as nonlinear problems such as non-constant coefficient Burger's equation. In the diffusion equation, the convergence rate is shown to be consistent with our analysis. That is noteworthy since we now have a theory guide for performing data collection and optimization. 

Based on the framework, we proposed a sensitivity analysis for the solution we obtained. The algorithm is based on the automatic differentiation mechanism and frees researchers from the tedious and error-prone process of deriving the gradients by hand. 

There are numerous research opportunities for solving inverse problems in differential equations based on neural networks and automatic differentiation
\begin{itemize}
    \item Tailored automatic differentiation software for engineering applications. Although the current software -- such as TensorFlow we have used for this article -- applies to a variety of problems, the performance is not tuned and optimized for engineering usage. Many subtle operators such as special treatment of boundary conditions will easily break the matrix calculation based software for deep learning nowadays. 
    \item Non-convex optimization for engineering problems. Although it is a well-known issue that problems involving neural networks may lead to non-convex optimization problems, based on authors' experience, such problems are less severe for the engineering problems considered compared to their counterparts in computer vision, speech recognition or natural language processing. Specific optimization techniques and theories may be developed for its own sake. 
    \item Neural network architectures for engineering problems. Currently, there are plenty of neural network architectures, and new research is going on. It would be interesting and fruitful to investigate the proper family of neural networks that are suitable for engineering problems, which have their unique properties~(such as known singularity at some locations). 
    \item Convergence theory for inverse problems using the neural networks. As we have seen in the analysis, it is possible to develop general theories about solving inverse problems using neural networks. Currently, our theory is case-dependent. A general theory is more attractive and beneficial for the development of the framework. 
    \item Uncertainty quantification. It is important to quantify the uncertainty under such framework since errors are inevitable in practical problems. We have proposed an tentative method to quantify the ``sensitivity'' using the automatic differentiation. The sensitivity region can be computed nearly for free thanks to the framework we have used: automatic differentiation is naturally embedded. However, a much more general uncertainty quantification approach shall be developed, and desirably leverages the automatic differentiation mechanism.
\end{itemize}

In conclusion, we believe the new framework will potentially bring useful tools for analysis and calibration of inverse problems in differential equations. 


\newpage

    \bibliographystyle{mybib}
\bibliography{sample.bib}

\end{document}